\documentclass[11pt,a4paper]{article}
\usepackage{indentfirst}
\usepackage{amsfonts}
\usepackage{amssymb}
\usepackage{mathrsfs}
\usepackage{amsmath}
\usepackage{amsthm}
\usepackage{enumerate}
\usepackage{cite}
\usepackage{tikz}
\usepackage{mathrsfs}
\allowdisplaybreaks
\usepackage{geometry}
\usepackage{ifpdf}
\ifpdf
  \usepackage[colorlinks=true,linkcolor=blue,citecolor=red, final,backref=page,hyperindex]{hyperref}
\else
  \usepackage[colorlinks,final,backref=page,hyperindex,hypertex]{hyperref}
\fi

\usepackage{txfonts}
\usepackage{indentfirst}
\usepackage{latexsym}

\usepackage[all]{xy}

\def\<{\langle}
\def\>{\rangle}

\def\c{\cdot}

\def\l{\lambda}

\newtheorem{thm}{Theorem}[section]
\newtheorem{lem}[thm]{Lemma}
\newtheorem{cor}[thm]{Corollary}
\newtheorem{pro}[thm]{Proposition}
\newtheorem{ex}[thm]{Example}

\theoremstyle{definition}
\newtheorem{defi}{Definition}[section]

\theoremstyle{remark}
\newtheorem{rmk}{Remark}[section]
\begin{document}
\title{\bf  On L-dendriform conformal 
algebras}
\author{\bf A. Hajjaji, L. Yuan}
\author{{  Atef Hajjaji$^{1}$
    \footnote {E-mail:  atef.hajjaji@uha.fr} \ and 
  Lamei Yuan$^{2}$
 \footnote { E-mail: lmyuan@hit.edu.cn (Corresponding author) }
}\\
{\small 1. IRIMAS - Département de Mathématiques, 18, rue des frères Lumière,
F-68093 Mulhouse, France. } \\
{\small 2.  School of Mathematics, Harbin Institute of Technology, Harbin 150001, China }}
\date{}
\maketitle
\begin{abstract}
In this paper, we introduce the concept of L-dendriform conformal algebras, which arise naturally from the study of $\mathcal{O}$-operators on left-symmetric conformal algebras and solutions to the conformal $S$-equation. These algebras extend the classical notions of dendriform and left-symmetric conformal algebras, providing a unified algebraic framework for understanding compatible structures in conformal algebra theory. We establish fundamental properties of L-dendriform conformal algebras, explore their relationships with $\mathcal{O}$ -operators, Rota-Baxter operators, and Nijenhuis operators, and demonstrate their connections to dendriform and quadri conformal algebras. Additionally, we investigate compatible $\mathcal{O}$-operators and their induced compatible L-dendriform conformal algebra structures. Our results generalize and unify several existing algebraic structures in conformal algebra theory, offering new insights into their interplay and applications. 
\end{abstract}



\numberwithin{equation}{section}

\tableofcontents

\section{Introduction} 
Lie conformal algebras, introduced by Kac \cite{KacVertex1996,KacFormal1999}, provide an algebraic framework for studying the singular part of operator product expansions in two-dimensional conformal field theory. They  
have deep connections with vertex algebras \cite{KacVertex1996} and numerous applications, for example in the Hamiltonian formalism of nonlinear evolution equations (see \cite{Dor} and the references therein, as well as \cite{GD1,ZELMANOV}). Considerable progress has also been made in the structure theory \cite{DK}, representation theory \cite{ChengConformal1997}, and cohomology theory \cite{BakalovCohomology1999} of finite Lie conformal algebras. Furthermore, in order to study whether compatible left-symmetric algebra structures exist on formal distribution Lie algebras, the notion of a left-symmetric conformal algebra was introduced in \cite{HL1}, which in turn provides a method for constructing vertex algebras. The notion of a finite left-symmetric conformal bialgebra which is free as
a $ \mathbb{C}[\partial] $-module was introduced and studied in \cite{HL2}.

As conformal analogues of left-symmetric algebras, left-symmetric conformal algebras have attracted significant attention due to their close relationship with Lie conformal algebras and their applications in integrable systems. They play a crucial role in understanding the algebraic structures underlying vertex algebras and classical Yang-Baxter equations \cite{Ku, BaiUnified2007,HongConformal2020}. The study of $\mathcal{O}$-operators on these algebras, as developed in \cite{BaiOOperators2012}, has revealed deep connections with Rota-Baxter operators and the conformal $S$-equation \cite{HongConformal2020}. Recently, many researchers have explored the cohomology and deformation theory of various operator algebras and conformal algebras, see \cite{YUAN1,YUAN2,ASIF,CHTIOUI,CHTIOUI2,CHTIOUI3,CHTIOUI4,SAMI}.

Dendriform algebras, introduced by Loday \cite{LodayDialgebras2002}, decompose the multiplication of an associative algebra into two operations, providing a combinatorial framework for algebraic structures. Their conformal analogues, dendriform conformal algebras, were studied in the context of Rota-Baxter operators and vertex algebras \cite{HB,YUAN1}. Quadri-algebras were introduced by Aguiar and Loday in \cite{AL} as a natural extension of Loday's work on dendriform algebras. Quadri conformal algebras, which serve as the conformal analogue of quadri-algebras and generalize dendriform conformal algebras, were introduced to study more complex algebraic structures \cite{AW}.


In this paper, we introduce {L-dendriform conformal algebras}, which serve as the conformal analogues of L-dendriform algebras introduced in \cite{BAILIU} as the underlying algebraic structure of a pseudo-Hessian pre-Lie algebra or a pseudo-Hessian Lie algebra. These structures emerge naturally from $\mathcal{O}$-operators on left-symmetric conformal algebras and provide solutions to the conformal $S$-equation. Our work is motivated by several important considerations:

\begin{itemize}
    \item L-dendriform conformal algebras serve as the underlying algebraic structures behind $\mathcal{O}$-operators of left-symmetric conformal algebras and the associated conformal $S$-equation. They provide explicit solutions to the conformal $S$-equation within certain left-symmetric conformal algebras.

    \item These algebras fit into a broader framework as Lie conformal analogues of dendriform conformal algebras, extending the classical relationship between Lie algebras and left-symmetric algebras to the conformal setting.

    \item There exist deep connections between $\mathcal{O}$-operators, Rota-Baxter operators, and Nijenhuis operators on left-symmetric conformal algebras. In particular, a Nijenhuis operator ``connects" two $\mathcal{O}$-operators in a way that any linear combination remains an $\mathcal{O}$-operator, leading naturally to compatible L-dendriform conformal algebra structures.
\end{itemize}

The paper is organized as follows: In Section \ref{s2}, we review left-symmetric conformal algebras, their modules, and $\mathcal{O}$-operators. In Section \ref{s3}, we introduce L-dendriform conformal algebras and studies their properties, including their relationships with dendriform and quadri conformal algebras. In Section \ref{s4}, we explore compatible $\mathcal{O}$-operators and their connections to compatible L-dendriform conformal algebras.

Our results reveal rich algebraic structures in the conformal setting and provide new insights into the interplay between different types of operators on left-symmetric conformal algebras. The constructions and theorems presented here extend known results from ordinary algebras to the conformal case and open new directions for further research in conformal algebra theory.

Throughout this paper, we denote by $\mathbb{C}$ the field of complex numbers. All tensors
over $\mathbb{C}$ are denoted by $\otimes$. We denote the identity map by $\mathrm{id}$. Moreover, if $A$ is a vector space, the space of polynomials
of $\l$ with coefficients in $A$ is denoted by $A[\lambda]$.
\section{Left-symmetric conformal algebras}\label{s2}
In this section, we begin by recalling some basic concepts related to Lie conformal algebras and left-symmetric conformal algebras. We then introduce the notion of $\mathcal{O}$-operators for left-symmetric conformal algebras in order to study the conformal $S$-equation. Furthermore, we show that an $\mathcal{O}$-operator of a left-symmetric conformal algebra associated with a bimodule gives rise to a symmetric solution of the conformal $S$-equation within a semi-direct product left-symmetric conformal algebra.

\subsection{Preliminaries and basic results}
In this subsection, we recall some definitions, notations and results about conformal
algebras. These materials can be found in \cite{DK,HL1,HL2,HongConformal2020,KacVertex1996,Li}. 
\begin{defi}
	A  conformal algebra  $A$ is a $ \mathbb{C}[\partial] $-module endowed with a $ \mathbb{C} $-bilinear map $ A\times A\rightarrow A[\lambda] $, denoted by $ a\times b\mapsto a \circ_{\lambda}b $
 satisfying the conformal sesquilinearity:
	\begin{eqnarray}\label{eq:conformal seq}
		(\partial a)\circ _{\lambda}b=-\lambda a \circ_{\lambda}b,\quad a\circ _{\lambda}(\partial b)=(\partial+\lambda)a\circ _{\lambda}b.
	\end{eqnarray}
 It will be denoted by $(A,\circ_\lambda)$. 
\end{defi}

A  conformal algebra is called 
 {\bf finite} if it is finitely generated as a $ \mathbb{C}[\partial] $-module. The
{\bf rank} of a conformal algebra $A$ is its rank as a $ \mathbb{C}[\partial] $-module.  

\begin{defi}
A  Lie conformal algebra $A$ is a $ \mathbb{C}[\partial] $-module endowed with a $ \mathbb{C} $-bilinear map $ A\times A\rightarrow A[\lambda] $, denoted by $ a\times b\mapsto [a_{\lambda}b] $ such that $A$ is a conformal algebra and 
\begin{eqnarray}
\label{eq:Lie conformal2}  {[a_{\lambda}b]}&=&-[b_{-\lambda-\partial}a],\\
 \label{eq:Lie conformal3} {[a_{\lambda}[b_{\mu}c]]}&=&[[a_{\lambda}b]_{\lambda+\mu}c]+[b_{\mu}[a_{\lambda}c]],\quad\forall~ a,b,c\in A.
\end{eqnarray} 
\end{defi}

\begin{defi}{\rm (\kern-3pt\cite{DK})
A module $M$ over a Lie conformal algebra ${A}$ is a $\mathbb{C}[\partial]$-module with a $\lambda$-action $a_\lambda m$, which is a $\mathbb{C}$-bilinear map ${A}\otimes M \rightarrow M[\lambda]=\mathbb{C}[\lambda]\otimes M$, $a\otimes m \mapsto a_\lambda m$, such that ($a,b\in {A}$ and $m\in M$)
\begin{align*}
&(\partial a)_\lambda m=-\lambda a_\lambda m,\quad a_\lambda(\partial m)=(\partial+\lambda)a_\lambda m,\\
&a_\lambda(b_\mu m)-b_\mu(a_\lambda m)=[a_\lambda b]_{\lambda+\mu}m.
\end{align*}}
\end{defi}

An $A$-module is called {\bf finite} if it is finitely generated as a $\mathbb{C} [\partial]$-module.

\begin{defi}
    Let $U$ and $V$ be two $\mathbb{C} [\partial]$-modules. A conformal linear map from
$U$ to $V$ is a $\mathbb{C}$-linear map $\phi:U\rightarrow V[\l]$, denoted by $\phi_{\l}:U\rightarrow V$, such that $[\partial,\phi_{\l}]=-\l\phi_{\l}$. Denote the $\mathbb{C}$-vector space of all such maps by Chom$(U, V)$. It has
a canonical structure of a $\mathbb{C} [\partial]$-module: 
\begin{equation}
    (\partial \phi)_{\l}=-\l \phi_{\l}.
\end{equation}
Define the conformal dual of a $\mathbb{C} [\partial]$-module $U$ as $U^{*c} = Chom(U,\mathbb{C} )$, where $\mathbb{C}$ is
viewed as the trivial $\mathbb{C} [\partial]$-module, that is,
\begin{equation*}
  U^{*c}=\{\psi:U\rightarrow  \mathbb{C}[\l] \;| \;\psi \text{ is a $\mathbb{C}$-linear and } \psi_{\l}(\partial b)=\l \psi_{\l}b,~~\forall~ b\in U\}
\end{equation*}
\end{defi}

Let $U$ and $V$ be finite modules over a Lie conformal algebra $A$. Then, the $\mathbb{C} [\partial]$-module Chom$(U, V)$ has an $A$-module structure defined by:
\begin{equation}
    (a_{\l}\phi)_{\mu}u=a_{\l}(\phi_{\mu-\l}u)-\phi_{\mu-\l}(a_{\l}u),
\end{equation}
for $a\in A$, $\phi\in Chom(U,V)$, $u \in U$. Hence, one special case is the contragradient
conformal $A$-module $U^{*c}$, where $\mathbb{C}$ is viewed as the trivial $A$-module and $\mathbb{C} [\partial]$-module. In particular, when $U = V$, set $Cend(V) = Chom(V, V)$. The $\mathbb{C} [\partial]$-module
$Cend(V)$ has a canonical structure of an associative conformal algebra defined by:
\begin{equation}
(a_{\l}b)_{\mu}v=a_{\l}(b_{\mu-\l}v),\quad a,b\in Cend(V),\;v\in V.
\end{equation}
Therefore, $gc(V) := Chom(V, V )$ has a Lie conformal algebra structure defined by:
\begin{equation}
    [a_{\l}b]_\mu v=a_{\l}(b_{\mu-\l}v)-b_{\mu-\l}(a_{\l}v),\quad a,b\in gc(V),\;v\in V.
\end{equation}
Here $gc(V )$ is called the {\bf general Lie conformal algebra} of $V$.\\

The tensor product $U \otimes V$ can be naturally endowed with an $A$-module structure
as follows:
\begin{equation*}
    \partial (u\otimes v)=\partial u\otimes  v+u\otimes \partial v,\quad a_{\l}(u\otimes v)=a_{\l}u\otimes v+u\otimes a_{\l}v,
\end{equation*}
where  $u\in U,\; v\in V$ and $a \in A$.
\begin{pro} (\cite{BKL1}) \label{BKL} Let $A$ be a Lie conformal algebra, and  
    let $U$ and $V$ be two $A$-modules. Suppose that $U$ is a $\mathbb{C} [\partial]$-module
of finite rank. Then, $U^{*c} \otimes V \cong Chom(U, V)$ as $A$-modules with the identification
$( f \otimes v)_{\l}(u) = f_{\l+\partial} (u)v$, where $f \in U^{*c}, u \in U$ and $v \in V$.
\end{pro}

\begin{defi}
A  left-symmetric conformal algebra is a $ \mathbb{C}[\partial] $-module $A$ endowed with a $ \mathbb{C} $-bilinear map $$ A\times A\rightarrow A[\lambda], ~~a\times b\mapsto a{\circ}_{\lambda}b, $$ such that $(A,\circ_\lambda)$ is a conformal algebra and the following identity hold for all $a,b,c\in A$:
\begin{equation}
\label{eq:LS conformal}   (a\circ _{\lambda}b)\circ_{\l+\mu}c-a\circ_{\l}(b\circ_{\mu}c)=(b\circ_{\mu}a)\circ_{\lambda+\mu}c-b\circ_{\mu}(a\circ_{\lambda}c).
\end{equation}
\end{defi}

\begin{ex}
    $A=\mathbb{C}[\partial]a$ is a free $\mathbb{C}[\partial]$-module of rank one. Define the $\l$-product on $A$ as follows: 
    \begin{equation*}
        a\circ_{\l}a=(\l+\partial +c)a,\quad c\in \mathbb{C}.
    \end{equation*}
    Then $A$ is a left-symmetric conformal algebra.
\end{ex}
\begin{ex}
Let $(A, \c)$ be a left-symmetric algebra. Then current left-symmetric conformal
algebra associated to $A$ is defined by:
\begin{equation*}
    Cur(A)=\mathbb{C}[\partial]\otimes A,\quad a\circ_{\l}b=a\c b,\;a,b\in A.
\end{equation*}
\end{ex}
\begin{pro}\label{LEFTLIE}
    If $A$ is a left-symmetric conformal algebra, then the $\l$-bracket
\begin{equation}\label{conformalcommutator}
  [a_{\l}b]=a\circ_{\l}b-b\circ_{-\partial-\l}a,\quad \forall a,b\in A,  
\end{equation}
defines a Lie conformal algebra $\mathfrak{g}(A)$, which is called {\bf the sub-adjacent Lie conformal algebra}
of $A$. In this case, $A$ is also called a {\bf compatible left-symmetric conformal algebra structure} on
the Lie conformal algebra $\mathfrak{g}(A)$.
\end{pro}
\begin{pro}
     Let $\mathfrak{g}$ be a $ \mathbb{C}[\partial] $-module  with a $\l$-binary operation $\circ_{\l}$. Then $(\mathfrak{g}, \circ_{\l})$ is a left-symmetric conformal 
algebra if and only if $(\mathfrak{g}, [\c_{\l}\c ])$ defined by Eq. \eqref{conformalcommutator} is a Lie conformal algebra.
\end{pro}
\begin{defi}
A module $V $ over a left-symmetric conformal algebra $A$ is a $\mathbb{C} [\partial]$-module with two $\mathbb{C}$-bilinear maps $A\times V\rightarrow V[\l],a\times v\mapsto a_{\l}v$ and $V\times A\rightarrow V[\l], v\times a\mapsto v_{\l}a$ such that  
	\begin{align}
		\label{eq:LS module1}  (\partial a)_{\lambda}v=-\lambda a_{\lambda}v,\quad \quad a_\lambda (\partial v)=(\partial+\lambda) a_\lambda v, \\
        \label{eq:LS module} 
        (\partial v)_{\lambda}a=-\lambda v_{\lambda}a, \quad \quad  v_\lambda (\partial a)=(\partial+\lambda) v_\lambda a,  \\
		\label{eq:LS module2}  (a\circ_{\lambda}b)_{\l+\mu} v-a_{\l}(b_{\mu}v)=(b\circ_{\mu}a)_{\lambda+\mu}v-b_{\mu}(a_{\l}v),\\
 \label{eq:LS module3}  
 (a_{\l}v)_{\l+\mu}b-a_{\l}(v_{\mu}b=(v_{\mu}a)_{\l+\mu}b-v_{\mu}(a\circ_{\l}b),
  	\end{align}
   hold for all $a, b \in A$ and $v \in V$.
\end{defi}
An $A$-module is called {\bf finite} if it is finitely generated as a $\mathbb{C} [\partial]$-module.
\begin{ex}
    Let $A$ be a left-symmetric algebra and $V$ be an $A$-module. Then $\mathbb{C} [\partial]V$ has a 
 natural bimodule structure over the current left-symmetric
conformal algebra $\mathbb{C} [\partial]A$, defined by:
\begin{equation*}
    a_{\l}v=a\c v,\quad v_{\l}a=v\c a, ~~ \forall ~ a\in A,~ v\in V,
\end{equation*}
where $a\c v$ and $v\c a$ denote left and right $A$-module actions, respectively, on $V$.
\end{ex}

\begin{rmk} Let $A$ be a left-symmetric conformal algebra and let $V$ be a finite $\mathbb{C}[\partial]$-module. Define $a_{\l}v = l_{A}(a)_{\l}v$ and $v_{\l}a = r_{A}(a)_{-\l-\partial}v$. It is straightforward to verify that a module structure of $V$ over the left-symmetric conformal algebra $A$ is equivalent to a pair ${l_A, r_A}$ of $\mathbb{C}[\partial]$-module homomorphisms $l_A, r_A : A \rightarrow \mathrm{Cend}(V)$ satisfying the following conditions for all $a, b \in A$ and $v \in V$: 
\begin{align}
    l_{A}(a\circ_{\l}b)_{\l+\mu}v-l_{A}(a)_{\l}(l_{A}(b)_{\mu}v)&=l_{A}(b\circ_{\mu}a)_{\l+\mu}v-l_{A}(b)_{\mu}(l_{A}(a)_{\l}v),\\
    r_{A}(b)_{-\partial-\l-\mu}(l_{A}(a)_{\l}v)-l_{A}(a)_{\l}(r_{A}(b)_{-\partial-\mu}v)&=  r_{A}(b)_{-\partial-\l-\mu}(r_{A}(a)_{\l}v)-r_{A}(a\circ_{\l}b)_{-\partial-\mu}v.
\end{align}
We denote this module by $(V,l_A, r_A)$.
\end{rmk}

Throughout this paper, we mainly consider $\mathbb{C}[\partial]$-modules that are finitely generated. For convenience, we will sometimes use the term \emph{representations of conformal algebras} to refer to such modules.

\begin{pro} (\kern-3pt\cite{HL2})\label{semidirect}
Let $A$ be a left-symmetric conformal algebra and $(V,l_A,r_A)$
be a representation of $A$ on $V$. Then, the $\mathbb{C} [\partial]$-module $A\oplus V$ is a left-symmetric conformal algebra
with the following $\l$-product: 
\begin{equation}\label{eqsemidirect}
    (a+u)\circ_{\l}(b+v)=a\circ_{\l}b+l_{A}(a)_{\l}v+r_{A}(b)_{-\partial-\l}u,\quad a,b\in A,\;u,v\in V.
\end{equation}
Denote it by $A \oplus_0 V$ , which is called the {\bf semi-direct product} of $A$ and $V$.
\end{pro}
\begin{pro} Let $A$ be a left-symmetric conformal algebra and
    $(V,l_A,r_A)$ a representation of $A$ on $V$. Then,
    \begin{enumerate}
        \item[(1)] $l_A : A \rightarrow gc(V)$ is a representation of the sub-adjacent Lie conformal algebra $\mathfrak{g}(A)$.
        \item[(2)] $\rho=l_A-r_A$ is a representation of the Lie conformal algebra $\mathfrak{g}(A)$.
    \end{enumerate}
\end{pro}
\begin{cor}
    Let $A$ be a  left-symmetric conformal algebra. Define two $\mathbb{C} [\partial]$-module homomorphisms $L_A$ and $R_A$ from $A$ to $Cend(A)$ by $L_{A}(a)_{\l}b = a_{\l}b$ and
$R_{A}(a)_{\l}b = b_{-\partial-\l}a$ for any $a, b \in A$. Then, $(A,L_A,R_A)$ is a representation of $A$ on $V$. Moreover, $L_A : A \rightarrow gc(A)$ and $\rho = L_A - R_A$ are
two representations of the Lie conformal algebra $\mathfrak{g}(A)$. 
\end{cor} 
\begin{pro}\label{repcoadj}
Let $A$ be a left-symmetric conformal algebra and $(V,l_A, r_A)$
be a representation of $A$ on $V$. Let $l^{*}_A$ and $r^{*}_A$
 be two $\mathbb{C} [\partial]$-module homomorphisms from $A$ to
$Cend(V^{*c})$ given by
\begin{equation}
  \big(l^{*}_A(a)_{\l}f\big)_{\mu}u=- f_{\mu-\l}\big(l_{A}(a)_{\l}u\big),\quad \big(r^{*}_A(a)_{\l}f\big)_{\mu}u=- f_{\mu-\l}\big(r_{A}(a)_{\l}u\big),  ~~\forall~ a\in A, ~ f\in V^{*c}, u\in V.  
\end{equation}
Then $(V^{*c},l^{*}_A-r^{*}_A,-r^{*}_A)$ is a representation of $A$ on $V^{*c}$.
\end{pro}

Next, let us introduce some notations about conformal bilinear form defined
in \cite{Li}.

Let $A$ be a $\mathbb{C} [\partial]$-module.    A {\bf conformal bilinear form} on $A$ is a $\mathbb{C}$-bilinear map $B_{\l}:A\otimes A\rightarrow \mathbb{C}[\l]$ satisfying
     \begin{equation}
       B_{\l}(\partial a,b)=-\l B_{\l}(a,b)=-B_{\l}(a,\partial b),\quad \forall a,b\in A.  
     \end{equation}
     If $$B_{\l}(a,b)=B_{-\l}(b,a), ~\forall~ a,b\in A,$$ we say this conformal bilinear form is {\bf symmetric}. 
     
     Suppose that $A$ is a free $\mathbb{C} [\partial]$-module of finite rank, and let $B_{\l}$ be a conformal bilinear
form on $A$. If the $\mathbb{C} [\partial]$-module homomorphism $T : A \rightarrow A^{*c}$, defined by $a \mapsto T_a $ with 
$$(T_a)_{\l}b =B_{\l}(a,b), ~\forall~b  \in A,$$ is an isomorphism, then we call the bilinear form {\bf nondegenerate}.

\begin{defi}
    Let $A$ be a left-symmetric conformal algebra. A conformal bilinear form $B_{\l} :
A \otimes A \rightarrow \mathbb{C} [\partial]$ is called a $2$-cocycle of $A$ if $B_{\l}$ satisfies the following property:
\begin{equation}
    B_{\l+\mu}(a\circ_{\l}b,c)-B_{\l}(a,b\circ_{\mu}c)=B_{\l+\mu}(b\circ_{\mu}a,c)-B_{\mu}(b,a\circ_{\l}c),
\end{equation}
for any $a, b, c \in A$.
\end{defi}
\begin{defi}(\kern-3pt\cite{HL2}) 
    Let $A$ be a
left-symmetric conformal algebra, and let $r=\sum_{i}x_i\otimes y_i\in A\otimes A$. Denote $\partial^{\otimes^3}=\partial \otimes 1\otimes 1+1\otimes \partial \otimes 1+1\otimes 1\otimes \partial$. The following equation
\begin{align*}
    &\{\{r,r\}\}:=\sum_{i,j}(y_{{j}_{\mu}}x_i\otimes x_j\otimes y_i)|_{\mu=1\otimes \partial \otimes 1}-\sum_{i,j}(x_j\otimes y_{{j}_{\mu}}x_i\otimes y_i)|_{\mu= \partial \otimes 1\otimes 1}\\
    &\quad \quad \quad -\sum_{i,j}(x_i\otimes x_j\otimes [y_{{i}_{\mu}}y_j])|_{\mu= \partial \otimes 1\otimes 1}=0\quad \text{mod }(\partial^{\otimes^3})\quad \text{in }A\otimes A \otimes A 
\end{align*}
is called the conformal $S$-equation in $A$.
\end{defi}
\subsection{$\mathcal{O}$-operators of left-symmetric conformal algebras and the conformal $S$-equation}

We introduce the notion of $\mathcal{O}$-operators for left-symmetric conformal algebras to interpret the conformal $S$-equation. In particular, an $\mathcal{O}$-operator of a left-symmetric conformal algebra  associated to a bimodule yields a symmetric solution of the conformal $S$-equation in a semi-direct product left-symmetric conformal algebra.

\begin{defi}
    Let $A$ be a
left-symmetric conformal algebra and $(V, l_A, r_A)$ be a bimodule of $A$. A $\mathbb{C} [\partial]$-module homomorphism $T:V\to A$ is called an $\mathcal{O}$-operator associated with $(V, l_A, r_A)$  if T satisfies
\begin{equation}
    T(u)\circ_{\l}T(v)=T\Big(l_A(T(u))_{\l}(v)+r_A(T(v))_{-\partial-\l}(u)\Big), \quad \forall u,v\in V. 
\end{equation}
In particular, an $\mathcal{O}$-operator $\mathcal{R}:A \to A$ associated with the adjoint bimodule $(A, L_A, R_A)$ is called a Rota-Baxter operator (of weight zero) on $A$; that is, $\mathcal{R}$ is a $\mathbb{C} [\partial]$-module homomorphism satisfying
\begin{equation}
    \mathcal{R}(a)\circ_{\l}\mathcal{R}(b)=\mathcal{R}\Big(\mathcal{R}(a)\circ_{\l}b+a\circ_{\l}\mathcal{R}(b)\Big), \quad \forall a,b\in A. 
\end{equation}
\end{defi}

Let $V$ be a bimodule over a left-symmetric conformal algebra $A$. By Proposition \ref{semidirect}, $A \oplus_0 V$ is a left-symmetric conformal algebra with respect to \eqref{eqsemidirect}. Suppose $T: V \rightarrow A$ is a $\mathbb{C}[\partial]$-module homomorphism, and  define the graph of $T$ as
\begin{equation}\label{eqgraph}
\operatorname{Gr}(T) = \{(T(v), v) \mid v \in V\}.
\end{equation}

\begin{pro} A $\mathbb{C}[\partial]$-module homomorphism 
    $T: V \rightarrow A$ is an $\mathcal{O}$-operator if and only if its graph 
$\operatorname{Gr}(T)$ is a subalgebra of $A \oplus_0 V$.
\end{pro} 
\begin{proof}
For any $(T(u), u), (T(v), v) \in \operatorname{Gr}(T)$, we have
\[
(T(u), u)\circ_\lambda (T(v), v) = \big(T(u)\circ_\lambda T(v), l_A(T(u))_\lambda v+ r_A(T(v))_{-\partial-\lambda}u\big).
\]
Hence, $T$ is an $\mathcal{O}$-operator if and only if 
$
\big(T(u)\circ_\lambda T(v), l_A(T(u))_\lambda v+ r_A(T(v))_{-\partial-\lambda}u\big) \in \operatorname{Gr}(T)[\lambda].$    \end{proof}

Note that $\operatorname{Gr}(T)$ and $V$ are isomorphic as $A$-bimodules by the identification $(T(u), u) \cong u$. Therefore, if $T$ is an $\mathcal{O}$-operator, i.e., if $\operatorname{Gr}(T)$ is a left-symmetric conformal subalgebra of $A \oplus_0 V$, then $V$ also inherits the structure of a left-symmetric conformal algebra.

Given an arbitrary $\mathbb{C}[\partial]$-module homomorphism $T: V \rightarrow A$, we define a lift of $T$, denoted by $\hat{T}$, as an endomorphism on $A \oplus V$ by
\[
\hat{T}(a, u) := (T(u), 0), \quad \text{for all } a \in A, u \in V.
\]

\begin{pro}\label{caraco-operator} A $\mathbb{C}[\partial]$-module homomorphism
$T: V \rightarrow A$ is an $\mathcal{O}$-operator if and only if  its lift $\hat{T}$ is a Rota-Baxter operator (of weight 0) on $A \oplus_0 V$.
   \end{pro}
\begin{proof} It is clear that 
$\hat{T}$ is  a $\mathbb{C}[\partial]$-module homomorphism. For any $(a, u), (b, v) \in A \oplus V$, we have
\begin{equation}\label{caraopera1}
\hat{T}(a, u)\circ_\lambda \hat{T}(b, v) = (T(u), 0)\circ_\lambda (T(v), 0) = (T(u)\circ_\lambda T(v), 0), 
\end{equation}
and
\begin{align}
\hat{T}\big( \hat{T}(a, u)\circ_\lambda (b, v) + (a, u)\circ_\lambda \hat{T}(b, v) \big)
&= \hat{T} \big( (T(u), 0)\circ_\lambda (b, v) + (a, u)\circ_\lambda (T(v), 0) \big) \notag \\
&= \hat{T} \big( (T(u)\circ_\lambda b, l_A(T(u))_\lambda v) + (a\circ_\lambda T(v), r_A(T(v))_{-\partial-\lambda}u \big) \notag \\
&= \big( T(l_A(T(u))_\lambda v+r_A(T(v))_{-\partial-\lambda}u), 0 \big). \label{caraopera2}
\end{align}
Combining \eqref{caraopera1} with \eqref{caraopera2}, we obtain the result.  
\end{proof}

Let \( A \) be a left-symmetric conformal algebra and \( r = \sum_i a_i \otimes b_i \in A \otimes A \). Define \( r^{21} = \sum_i b_i \otimes a_i \).  
We say \( r \) is \emph{skew-symmetric} if \( r = -r^{21} \), and \emph{symmetric} if \( r = r^{21} \). Let \( V \) be a free \( \mathbb{C} [\partial] \)-module of finite rank and let \( \{ v_i \}_{i=1}^m \) be a \( \mathbb{C} [\partial] \)-basis of \( V \).  
If $V$ is a bimodule of \( A \), then \( V \cong {V^{*c}}^{*c} \) via the \( \mathbb{C} [\partial] \)-module homomorphism \( v_i \mapsto v_i^{**} \), where  
\[v^{**}_{{i}_{\;\;-\lambda - \partial}}(y^*)= y^*_\lambda(v_i), \quad \text{for any } y^* \in V^{*c}. \]

Now suppose that \( A \) is a finite left-symmetric conformal algebra that is free as a \( \mathbb{C}[\partial] \)-module.  
Then, by the above discussion, we have the following isomorphisms of \( A \)-bimodules:  
\[ A \otimes A \cong {A^{*c}}^{*c} \otimes A \cong C\mathrm{hom}(A^{*c}, A). \]

For any \( a, b \in A \) and \( u, v \in A^{*c} \), we define 
\[
\{u, a\}_\lambda = u_\lambda(a), ~\text{and} ~ \{u \otimes v, a \otimes b\}_{(\lambda, \mu)} = \{u, a\}_\lambda \{v, b\}_\mu.
\]
Then, by Proposition \ref{BKL}, for any \( r \in A \otimes A \), we associate to it a conformal linear map  
\( T^r \in C\mathrm{hom}(A^{*c}, A) \) defined as follows:
\begin{equation}\label{eqpair}
  \{f, T^{r}_{-\mu - \partial}(g)\}_\lambda = \{g \otimes f, r\}_{(\mu, \lambda)}, \quad \forall ~ f, g \in A^{*c}.  
\end{equation}
 By Proposition \ref{repcoadj}, we have
 \begin{equation*}
  <a\circ_{\l} b,u>_\mu= <b,L_{A}^{*}(a)_{\l}>_{\mu-\l},\quad <b\circ_{\mu-\l} a,u>_\mu= <b,R_{A}^{*}(a)_{\l}>_{\mu-\l}, \; \forall a,b\in A, u\in A^{*c}.  
 \end{equation*}
Set \( r = \sum_i a_i \otimes b_i \in A \otimes A \). Then, from Eq. \eqref{eqpair} we obtain
\[
T^{r}_\lambda(u) = \sum_i \{u, a_i\}_{-\lambda - \partial} b_i, \quad \forall~ u \in A^{*c}.
\]

\begin{thm}\label{solofS-eq}
    Let $A$ be a free left-symmetric conformal algebra of finite rank
and let $r \in A \otimes A$ be symmetric. Then, $r$ is a solution of the conformal $S$-equation if and
only if $T^r \in C\mathrm{hom}(A^{*c}, A)$ corresponding to $r$ satisfies
\begin{equation}
    T^{r}_0(u)\circ_{\l}T^{r}_0(v)=T^{r}_0\Big((L_A^{*}-R_A^{*})(T^{r}_0(u))_{\l}(v)-R_A^{*}(T^{r}_0(v))_{-\partial-\l}(u)\Big), \quad \forall~~ u,v\in A^{*c}, 
\end{equation}
where \( T^{r}_0 = T^{r}_\lambda\; \big|_{\lambda=0} \).
\end{thm}
\begin{proof} For any \( a \in A \) and \( u \in A^{*c} \), 
define
\[
\langle a, u\rangle_\lambda = \{u, a\}_{-\lambda} = u_{-\lambda}(a).
\]
Clearly, we have
\begin{equation}\label{eq2.23}
\langle \partial a, u\rangle_\lambda = -\lambda \langle a, u\rangle_\lambda. 
\end{equation}
Moreover, define
\[
\langle a \otimes b \otimes c, u \otimes v \otimes w\rangle_{(\lambda, \nu, \theta)} = \langle a, u\rangle _\lambda \langle b, v\rangle_\nu \langle c, w\rangle_\theta,
\]
where \( a, b, c \in A \), and \( u, v, w \in A^{*c} \). Let \( r = \sum_i a_i \otimes b_i \in A \otimes A \) be symmetric. From the discussion above, the map \( T^r \in \operatorname{Chom}(A^{*c}, A) \) corresponding to \( r \) is given by:
\[
T^r_\lambda(u) = \sum_i \langle a_i, u\rangle _{\lambda + \partial} b_i, \quad u \in A^{*c}.
\]
Since \( r \) is symmetric, one has
\[
T^r_\lambda(v) =  \sum_i \langle b_i, v\rangle _{\lambda + \partial} a_i, \quad v \in A^{*c}.
\]

Now suppose that $r$ is a solution of the conformal $S$-equation, i.e., $r$ satisfies  
\begin{equation}\label{eq2.25}
\langle \{\{r, r\}\} \mod (\partial^{\otimes 3}), u \otimes v \otimes w\rangle_{(\lambda, \nu, \theta)} = 0, \quad \text{for any } u, v, w \in A^{*c}. 
\end{equation}
By Eq. \eqref{eq2.23}, the condition in Eq. \eqref{eq2.25} is equivalent to:
\begin{equation}
\big\langle\{\{r, r\}\}, u \otimes v \otimes w\big\rangle_{(\lambda, \nu, \theta)} = 0 \mod (\lambda + \nu + \theta). 
\end{equation}
By a direct computation, we obtain 
\begin{eqnarray*}
\Big\langle\sum_{i,j} b_j\circ_\mu a_i \otimes a_j \otimes b_i \big|_{\mu = 1 \otimes \partial \otimes 1}, u \otimes v \otimes w\Big\rangle_{(\lambda, \nu, \theta)}
&=&\sum_{i,j} \langle b_j\circ_{-\nu} a_i,u\rangle_{\lambda}\langle a_j,v\rangle _{\nu}\langle b_i,v\rangle_{\theta}\\
&=&\Big\langle \Big(\sum_{j} \langle a_j, v\rangle_\nu b_j\Big)\circ_{-\nu} \Big(\sum_{i} \langle b_i, w\rangle_\theta a_i\Big), u\Big\rangle_\lambda \\
&=& \Big< T^r_{\nu-\partial}(v)\circ_{-\nu} T^r_{\theta-\partial}(w), u\Big>_\lambda\\&=&\Big< T^r_{0}(v)\circ_{-\nu} T^r_{\l+\nu+\theta}(w), u\Big>_\lambda.    
\end{eqnarray*}
Similarly, we have
\[
\Big\langle\sum_{i,j} a_j \otimes b_j\circ_\mu a_i \otimes b_i \big|_{\mu = \partial \otimes 1 \otimes 1}, u \otimes v \otimes w\Big\rangle_{(\lambda, \nu, \theta)}
= -\langle T^r_{\lambda + \nu + \theta}(R_A^*(T^r_{0}(w))_{-\theta} v), u\rangle_\l,
\]
and
\begin{eqnarray*}
\Big\langle \sum_{i,j} a_i \otimes a_j \otimes [b_i{}_\mu b_j] \big|_{\mu = \partial \otimes 1 \otimes 1}, u \otimes v \otimes w\Big\rangle_{(\lambda, \nu, \theta)}
= -\langle T^r_{\lambda + \nu + \theta}(R_A^*(T^r_{0}(v))_{-\nu} w), u\rangle_\l
+\langle T^r_{\lambda + \nu + \theta}(L_A^*(T^r_0(v))_{-\nu} w), u\rangle_\l.
\end{eqnarray*}
Using Eq. \eqref{eq2.25} and the above discussion, the conformal $S$-equation is equivalent to:
\begin{align}\label{eq2.26}
\Big< T^r_{0}(v)\circ_{-\nu} T^r_{\l+\nu+\theta}(w)& +T^r_{\lambda + \nu + \theta}(R_A^*(T^r_{0}(w))_{-\theta} v)+T^r_{\lambda + \nu + \theta}(R_A^*(T^r_{0}(v))_{-\nu} w)\nonumber\\
&-T^r_{\lambda + \nu + \theta}(L_A^*(T^r_0(v))_{-\nu} w),u\Big>_{\l}= 0 \mod (\lambda + \nu + \theta).
\end{align}
This is equivalent to:
\begin{equation}
T^{r}_0(v)\circ_{-\nu}T^{r}_0(w)=T^{r}_0\Big((L_A^{*}-R_A^{*})(T^{r}_0(v))_{-\nu}(w)-R_A^{*}(T^{r}_0(w))_{-\partial+\nu}(v)\Big).
\end{equation}
Finally, the conclusion is obtained directly by replacing \( -\nu \) by \( \nu \).
\end{proof}

\begin{rmk} Let $A$ be a free left-symmetric conformal algebra of finite  rank.  
    By Theorem \ref{solofS-eq}, a symmetric solution of conformal S-equation in $A$ is
equivalent to the existence of $T \in Chom(A^{*c}, A)$, such that $T^r_{0}$ is an $\mathcal{O}$-operator associated with $\big(L_A^{*}-R_A^{*},-R_A^{*}\big)$.
\end{rmk}


For any vector space \( V \), let 
$\tau : V \otimes V \rightarrow V \otimes V$
be the {\bf flip map}, that is,
\[
\tau(x \otimes y) = y \otimes x, \quad \forall\, x, y \in V.
\]
Next, we study the $\mathcal{O}$-operators of left-symmetric conformal algebras in a more general setting. Let $(V, l_A, r_A)$ be a bimodule of a left-symmetric conformal algebra $A$.  By Proposition \ref{repcoadj}, $(V^{*c}, l_A^*-r_A^*, -r_A^*)$ is a bimodule of $A$. Suppose that $V$ is a $ \mathbb{C}[\partial] $-module of finite rank. By Proposition \ref{BKL}, $V^{*c} \otimes A \cong \operatorname{Chom}(V, A)$ as $\mathbb{C}[\partial]$-modules via the isomorphism  $\varphi$ defined by
\[
\varphi(f \otimes a)_\lambda v = f_{\lambda + \partial^A}(v) a, \quad \forall a \in A, v \in M, f \in V^{*c}.
\]
Moreover, due to the $\mathbb{C}[\partial]$-module structures on $V^{*c} \otimes A$, we also have: $V^{*c} \otimes A \cong A \otimes V^{*c}$
as $\mathbb{C}[\partial]$-modules. Then,  $$\operatorname{Chom}(V, A) \cong A \otimes V^{*c},$$ as $\mathbb{C}[\partial]$-modules. Therefore, for any 
$T \in \operatorname{Chom}(V, A)$, we can associate an element $r_T \in A \otimes V^{*c} \subset (A \ltimes_{l_A^*-r_A^*,-r_A^* } V^{*c}) \otimes (A \ltimes_{l_A^*-r_A^*,-r_A^*} V^{*c})$.

\medskip

\begin{thm}\label{solS-equa}
Let $A$ be a finite left-symmetric conformal algebra and $(V, l_A, r_A)$ be a finite bimodule of $A$. Suppose that $A$ and $V$ are free as $\mathbb{C}[\partial]$-modules. Let 
$T \in \operatorname{Chom}(V, A)$ and $r_T \in A \otimes V^{*c} \subset (A \ltimes_{l_A^*-r_A^*,-r_A^* } V^{*c}) \otimes (A \ltimes_{l_A^*-r_A^*,-r_A^*} V^{*c})$ be the element corresponding to $T$ under the above isomorphism. Then
\begin{align*}
 r = r_T + \tau r_T
\end{align*}
is a symmetric solution of the conformal $S$-equation in the left-symmetric conformal algebra  
$A \ltimes_{l_A^*-r_A^*,-r_A^*} V^{*c}$ if and only if $T_0 := T_\lambda|_{\lambda = 0}$ is an $\mathcal{O}$-operator associated with the bimodule $(V, l_A, r_A)$.
\end{thm}

\begin{proof}
 Let $\{e_1, \dots, e_n\}$ be a $\mathbb{C}[\partial]$-basis of $A$, $\{v_1, \dots, v_m\}$ be a $\mathbb{C}[\partial]$-basis of $V$, and $\{v_1^*, \dots, v_m^*\}$ be the dual  $\mathbb{C}[\partial]$-basis of $V^{*c}$. Suppose
\[
T_\lambda(v_i) = \sum_{j=1}^n a_{ij}(\lambda, \partial) e_j
, \quad \forall i = 1, \dots, m,
\]
with $a_{ij}(\lambda, \partial) \in \mathbb{C}[\lambda, \partial]$. Then
\[
r_T = \sum_{i=1}^m \sum_{j=1}^n a_{ij}(-1 \otimes \partial - \partial \otimes 1, \partial \otimes 1) \, e_j \otimes v_i^*.
\]
Hence,
\[
r = \sum_{i,j} a_{ij}(-\partial \otimes 1 - 1 \otimes \partial, \partial \otimes 1) \, e_j \otimes v_i^*
 + \sum_{i,j} a_{ij}(-\partial \otimes 1 - 1 \otimes \partial, 1 \otimes \partial) \, v_i^* \otimes e_j.
\]
Moreover, from the definitions of $l_A^*$ and $r_A^*$, we have:
\begin{equation}\label{defofl*r*}
 l_A^*(e_i)_\lambda v_j^* =- \sum_k v_j^*{}_{-\lambda - \partial}(l_A(e_i)_\lambda v_k)v_k^*, \quad r_A^*(e_i)_\lambda v_j^* = -\sum_k v_j^*{}_{-\lambda - \partial}(r_A(e_i)_\lambda v_k)v_k^*.   
\end{equation}
Then we get:
\begin{align*}
\{\{r, r\}\} = & 
\left(
\sum_{i,j,k,l} a_{ij}(0, -1 \otimes \partial \otimes 1) a_{kl}(0, -1 \otimes 1 \otimes \partial)\, e_l {\circ}_\mu e_j \otimes v_k^* \otimes v_i^* \right. \\
&\left. + \sum_{i,j,k,l} a_{ij}(0, 1 \otimes \partial \otimes 1) a_{kl}(0, -1 \otimes 1 \otimes \partial)\, v_k^* {\circ}_\mu e_j \otimes e_l \otimes v_i^* \right. \\
&\left. + \sum_{i,j,k,l} a_{ij}(0, -1 \otimes \partial \otimes 1) a_{kl}(0, 1 \otimes 1 \otimes \partial)\, e_l {\circ}_\mu v_i^*] \otimes v_k^* \otimes e_j
\right)\Big|_{\mu=1 \otimes \partial \otimes 1} \\
& - \left(
\sum_{i,j,k,l} a_{ij}(0, \partial \otimes 1 \otimes 1) a_{kl}(0, -1 \otimes 1 \otimes \partial)\, v^*_k \otimes e_l {\circ}_\mu v_i^* \otimes e_j \right. \\
&\left. + \sum_{i,j,k,l} a_{ij}(0, -\partial \otimes 1 \otimes 1) a_{kl}(0, -1 \otimes 1 \otimes \partial)\, v_k^* \otimes e_l {\circ}_\mu e_j \otimes v_i^* \right. \\
&\left. + \sum_{i,j,k,l} a_{ij}(0, -\partial \otimes 1 \otimes 1) a_{kl}(0, 1 \otimes 1 \otimes \partial)\, e_l \otimes v_k^* {\circ}_\mu e_j \otimes v^*_i
\right)\Big|_{\mu=\partial \otimes 1 \otimes 1} \\
& - \left(
 \sum_{i,j,k,l} a_{ij}(0, \partial \otimes 1 \otimes 1) a_{kl}(0, -1 \otimes \partial \otimes 1)\, v_i^* \otimes e_l \otimes e_j {\circ}_\mu v_k^* \right. \\
&\left. + \sum_{i,j,k,l} a_{ij}(0, -\partial \otimes 1 \otimes 1) a_{kl}(0, 1 \otimes \partial \otimes 1)\,  e_j\otimes v_k^* \otimes v_i^* {\circ}_\mu e_l \right. \\
&\left. + \sum_{i,j,k,l} a_{ij}(0, -\partial \otimes 1 \otimes 1) a_{kl}(0, -1 \otimes \partial \otimes 1)\, v_i^* \otimes v_k^* \otimes e_j {\circ}_\mu e_l
\right)\Big|_{\mu=\partial \otimes 1 \otimes 1} \\%
& + \left(
 \sum_{i,j,k,l} a_{ij}(0, \partial \otimes 1 \otimes 1) a_{kl}(0, -1 \otimes \partial \otimes 1)\, v_i^* \otimes e_l \otimes v_k^* {\circ}_{-\partial-\mu} e_j \right. \\
&\left. + \sum_{i,j,k,l} a_{ij}(0, -\partial \otimes 1 \otimes 1) a_{kl}(0, 1 \otimes \partial \otimes 1)\,  e_j\otimes v_k^* \otimes e_l {\circ}_{-\partial-\mu} v_i^* \right. \\
&\left. + \sum_{i,j,k,l} a_{ij}(0, -\partial \otimes 1 \otimes 1) a_{kl}(0, -1 \otimes \partial \otimes 1)\, v_i^* \otimes v_k^* \otimes e_l {\circ}_{-\partial-\mu} e_j
\right)\Big|_{\mu=\partial \otimes 1 \otimes 1}\mod (\partial^{\otimes 3}).
\end{align*}

Note that $T_0(v_i) = \sum_{j=1}^n a_{ij}(0, \partial) e_j,$ from which we obtain

\begin{align*}
\{\{r, r\}\} = \sum_{i,k} \bigg(&
T_0(v_k) {\circ}_\mu T_0(v_i) \otimes v_k^* \otimes v_i^* 
+ (-r_A^*)(T_0(v_i))_{-\mu - \partial} v_k^* \otimes T_0(v_k) \otimes v_i^* \\
&+ (l_A^*-r_A^*)(T_0(v_k))_\mu v_i^* \otimes v_k^* \otimes T_0(v_i)
\bigg)\Big|_{\mu = 1 \otimes \partial \otimes 1} \\
- \sum_{i,k} \bigg(&
T_0(v_k) \otimes -r_A^*(T_0(v_k))_{-\mu - \partial} v_k^* \otimes v_i^*
+ v_k^* \otimes T_0(v_k) {\circ}_\mu T_0(v_i) \otimes v_i^* \\
&+ v_k^* \otimes (l_A^*-r_A^*)(T_0(v_k))_{\mu} v_i^* \otimes T_0(v_i)
\bigg)\Big|_{\mu = \partial \otimes 1 \otimes 1} \\
- \sum_{i,k} \bigg(&
 v_i^* \otimes T_0(v_k) \otimes l_A^*(T_0(v_i))_{\mu} v_k^*- T_0(v_i) \otimes v_k^* \otimes l_A^*(T_0(v_k))_{-\partial-\mu} v_i^*\\
&+ v_i^* \otimes v_k^* \otimes T_0(v_i) {\circ}_\mu T_0(v_k)-v_i^* \otimes v_k^* \otimes T_0(v_k) {\circ}_{-\partial-\mu} T_0(v_i)
\bigg)\Big|_{\mu = \partial \otimes 1 \otimes 1}
\mod (\partial^{\otimes 3}).
\end{align*}

Since \( T_0 \) is a \( \mathbb{C}[\partial] \)-module homomorphism and by Eq. \eqref{defofl*r*}, we have:
\begin{align*}
&\sum_{i,k} T_0(v_i) \otimes v_k^* \otimes l_A^* (T_0(v_k))\,{}_{-\mu - \partial} v_i^* \Big|_{\mu = \partial \otimes 1 \otimes 1} \\
\equiv& \sum_{i,k} T_0(v_i) \otimes v_k^* \otimes l_A^*(T_0(v_k))\,{}_{1 \otimes \partial \otimes 1} v_i^* \mod (\partial^{\otimes 3}) \\
\equiv& -\sum_{i,j,k} T_0(v_i) \otimes v_k^* \otimes v_i^*\,_{\partial \otimes 1 \otimes 1} (l_A(T_0(v_k))\,{}_{1 \otimes \partial \otimes 1} v_j)\, v_j^* \mod (\partial^{\otimes 3}) \\
\equiv & -\sum_{i,j,k} T_0\left(v_i^*\, _{\partial}(l_A(T_0(v_k))\,{}_{1 \otimes \partial \otimes 1} v_j)\, v_i\right) \otimes v_k^* \otimes v_j^* \mod (\partial^{\otimes 3}) \\
\equiv &-\sum_{i,k} T_0\left(l_A(T_0(v_k))\,{}_\mu v_i \right) \otimes v_k^* \otimes v_i^* \Big|_{\mu = 1 \otimes \partial \otimes 1} \mod (\partial^{\otimes 3}).
\end{align*}

Similarly, we have:
\begin{align*}
&\{\{r, r\}\} \mod (\partial^{\otimes 3})\\
\equiv& \sum_{i,k} \Big(
 \big( T_0(v_k)\,{\circ}_\mu T_0(v_i) - T_0(l_A(T_0(v_k))\,{}_\mu v_i) - T_0(r_A(T_0(v_i))\,{}_{-\mu - \partial} v_k) \big) \otimes v_k^* \otimes v_i^* \Big|_{\mu = 1 \otimes \partial \otimes 1} \\
& -  \, v_k^* \otimes \big( T_0(v_k)\,{\circ}_\mu T_0(v_i) - T_0(l_A(T_0(v_k))\,{}_\mu v_i) - T_0(r_A(T_0(v_i))\,{}_{-\mu - \partial} v_k) \big)\otimes v_i^* \Big|_{\mu = \partial  \otimes 1 \otimes 1} \\
&-  \, v_k^* \otimes v_i^* \otimes \big( T_0(v_k)\,{\circ}_\mu T_0(v_i) - T_0(l_A(T_0(v_k))\,{}_\mu v_i) - T_0(r_A(T_0(v_i))\,{}_{-\mu - \partial} v_k) \big) \Big|_{\mu = \partial \otimes 1 \otimes 1}\\
&+  \, v_i^* \otimes v_k^* \otimes \big( T_0(v_k)\,{\circ}_{-\mu-\partial} T_0(v_i) - T_0(l_A(T_0(v_k))\,{}_{-\mu-\partial} v_i) - T_0(r_A(T_0(v_i))\,{}_{\mu} v_k) \big) \Big|_{\mu = \partial \otimes 1 \otimes 1}\Big).
\end{align*}
Therefore, \( r \) is a solution of the conformal $S$-equation in the left-symmetric conformal $A \ltimes_{l_A^*-r_A^*,-r_A^*} V^{*c}$ if and only if 
\[
T_0(v_k)\,{}_\mu T_0(v_i) = T_0(l_A(T_0(v_k))\,{}_\mu v_i) + T_0(r_A(T_0(v_i))\,{}_{-\mu - \partial} v_k), \quad \forall i, k \in \{1, \dots, m\}.
\]
This completes the proof.  
\end{proof}   

\section{L-dendriform conformal algebras }\label{s3}
In this section, we introduce the concept of an L-dendriform conformal algebra and investigate its fundamental properties in terms of $\mathcal{O}$-operators on left-symmetric conformal algebras. In particular, we examine their connections with left-symmetric conformal algebras, conformal $S$-equations, pseudo-Hessian structures, and other related conformal algebraic structures.
\subsection{The definition and some basic properties}
\begin{defi}
An L-dendriform conformal algebra $A$ is a $ \mathbb{C}[\partial] $-module endowed with two bilinear $\l$-operations $$ \vartriangleright_{\l}:A\times A\rightarrow A[\l]~~\text{and} ~~ \vartriangleleft_{\l}: A\times A\rightarrow A[\lambda], $$ satisfying the conformal sesquilinearity and the following axioms for all $a, b, c \in A$:
\begin{eqnarray}
\label{eq1:L-dendriform conformal}   
a\vartriangleright_{\lambda}(b\vartriangleright_{\mu}c)&=&(a\vartriangleright_{\l}b)\vartriangleright_{\l+\mu}c+(a\vartriangleleft_{\l}b)\vartriangleright_{\l+\mu}c+b\vartriangleright_{\mu}(a\vartriangleright_{\l}c)\nonumber\\
&&-\;(b \vartriangleright_{\mu}a)\vartriangleright_{\lambda+\mu}c-(b \vartriangleleft_{\mu}a)\vartriangleright_{\lambda+\mu}c,\\
\label{eq2:L-dendriform conformal}
a\vartriangleright_{\lambda}(b\vartriangleleft_{\mu}c)&=&(a\vartriangleright_{\l}b)\vartriangleleft_{\l+\mu}c+b\vartriangleleft_{\mu}(a\vartriangleright_{\l}c)+b\vartriangleleft_{\mu}(a\vartriangleleft_{\l}c)\nonumber\\
&&-\;(b \vartriangleleft_{\mu}a)\vartriangleleft_{\lambda+\mu}c.
\end{eqnarray}
\end{defi}

\begin{pro}\label{denleft}
  Let $(A, \vartriangleright_{\l}, \vartriangleleft_{\l})$  be an L-dendriform conformal algebra.
  \begin{enumerate}
      \item[(1)] The $\l$-operation $ \bullet_{\l}: A\times A\rightarrow A[\lambda] $ defined by
      \begin{equation}\label{eqbullet}
        a \bullet_{\l}b=a \vartriangleright_{\l}b+a \vartriangleleft_{\l}b,\quad\forall~ a,b\in A,
      \end{equation}
gives $A$ a left-symmetric conformal algebra structure. The algebra  $(A,\bullet_{\l})$ is called the {\bf associated horizontal left-symmetric conformal algebra} of $(A, \vartriangleright_{\l}, \vartriangleleft_{\l})$ and  we say that $(A, \vartriangleright_{\l}, \vartriangleleft_{\l})$ is  a {\bf 
compatible L-dendriform conformal algebra structure} on $(A,\bullet_{\l})$.
\item[(2)] The $\l$-operation $ \circ_{\l}: A\times A\rightarrow A[\lambda] $ defined by
      \begin{equation}\label{eqcirc}
        a \circ_{\l}b=a \vartriangleright_{\l}b-b \vartriangleleft_{-\partial-\l}a,\quad\forall~ a,b\in A,
      \end{equation}
gives $A$ the structure of 
 a left-symmetric conformal algebra. The algebra $(A,\circ_{\l})$ is called the {\bf associated vertical left-symmetric conformal algebra} of $(A, \vartriangleright_{\l}, \vartriangleleft_{\l})$ and we say that $(A, \vartriangleright_{\l}, \vartriangleleft_{\l})$ is a
{\bf compatible L-dendriform conformal algebra structure} on $(A,\circ_{\l})$.
\item[(3)] Both $(A,\bullet_{\l})$ and $(A,\circ_{\l})$ have the same {\bf sub-adjacent Lie conformal algebra $\mathfrak{g}(A)$}, whose $\l$-bracket is given by 
\begin{equation}
    [a_{\l}b]=a\vartriangleright_{\l}b+a \vartriangleleft_{\l}b-b \vartriangleright_{-\partial-\l}a-b \vartriangleleft_{-\partial-\l}a,\quad\forall~ a,b\in A.
\end{equation}
  \end{enumerate}
\end{pro}
\begin{proof}
  The conclusions can be obtained by straightforward computation.  
\end{proof}

\begin{rmk}\label{genass}
   Let $(A, \vartriangleright_{\l}, \vartriangleleft_{\l})$  be an L-dendriform conformal algebra. Then the defining relations \eqref{eq1:L-dendriform conformal} and \eqref{eq2:L-dendriform conformal} can be rewritten as (for any $a,b,c\in A$)
   \begin{eqnarray}
      \label{eq1:L-dendriform conformal1}   
a\vartriangleright_{\lambda}(b\vartriangleright_{\mu}c)-(a\bullet_{\l}b)\vartriangleright_{\l+\mu}c&=&b\vartriangleright_{\mu}(a\vartriangleright_{\l}c)-(b \bullet_{\mu}a)\vartriangleright_{\lambda+\mu}c,\\
\label{eq2:L-dendriform conformal2}
a\vartriangleright_{\lambda}(b\vartriangleleft_{\mu}c)-(a\vartriangleright_{\l}b)\vartriangleleft_{\l+\mu}c&=&b\vartriangleleft_{\mu}(a\bullet_{\l}c)-(b \vartriangleleft_{\mu}a)\vartriangleleft_{\lambda+\mu}c, 
   \end{eqnarray}
where $\bullet_{\lambda}$ is defined as in Eq.~\eqref{eqbullet}.  
Both sides of these equations can be viewed as certain kinds of \textbf{generalized associators}. In this sense, Eqs.~\eqref{eq1:L-dendriform conformal1} and \eqref{eq2:L-dendriform conformal2} express a form of \textbf{generalized conformal left-symmetry} for these generalized associators.
\end{rmk}

\begin{pro}\label{horver}
    Let $A$ be a $ \mathbb{C}[\partial] $-module endowed with two bilinear $\l$-operations
    
    $$ \vartriangleright_{\l}:A\times A\rightarrow A[\l]~~\text{and} ~~ \vartriangleleft_{\l}: A\times A\rightarrow A[\lambda]. $$
    \begin{enumerate}
        \item[(1)]  $(A, \vartriangleright_{\l}, \vartriangleleft_{\l})$  is an L-dendriform conformal algebra if and only if $(A,\bullet_{\l})$ defined by Eq. \eqref{eqbullet} is a left-symmetric conformal algebra and $(A,L_{\vartriangleright_{\l}},R_{\vartriangleleft_{\l}})$ is a bimodule.
\item[(2)] $(A, \vartriangleright_{\l}, \vartriangleleft_{\l})$ is an L-dendriform algebra if and only if $(A,\circ_{\l})$  defined by Eq. \eqref{eqcirc} is a left-symmetric conformal algebra and $(A,L_{\vartriangleright_{\l}},-L_{\vartriangleleft_{\l}})$ is a bimodule.
  \end{enumerate}
\end{pro}
\begin{proof}
The conclusions can be obtained by straightforward computation or a proof similar to that of Theorem \ref{thmdenonVver}.    
\end{proof}
\begin{cor}\label{corhorver}
    Let $(A, \vartriangleright_{\l}, \vartriangleleft_{\l})$ be an L-dendriform conformal algebra. Then $(A^{*c},L^{*}_{\vartriangleright_{\l}}-R^{*}_{\vartriangleleft_{\l}},-R^{*}_{\vartriangleleft_{\l}})$ is a bimodule of the associated horizontal left-symmetric conformal
algebra $(A,\bullet_{\l})$ and $(A^{*c},L^{*}_{\vartriangleright_{\l}}+L^{*}_{\vartriangleleft_{\l}},L^{*}_{\vartriangleleft_{\l}})$ is a bimodule of the associated vertical left-symmetric conformal algebra $(A,\circ_{\l})$.
\end{cor}
\begin{proof}
    It follows from Propositions \ref{repcoadj} and \ref{horver}.
\end{proof}

\begin{pro}
 Let $(A, \vartriangleright_{\l}, \vartriangleleft_{\l})$ be an L-dendriform conformal algebra. Define two new binary $\l$-operations $\vartriangleright^{t}_{\l}$ and $ \vartriangleleft^{t}_{\l}: A\times A\rightarrow A[\lambda] $ by
\begin{equation}\label{tranoperations}
 a\vartriangleright^{t}_{\l}b= a\vartriangleright_{\l}b, \quad \quad  a\vartriangleleft^{t}_{\l}b=-b\vartriangleleft_{-\partial-\l}a,\quad \forall~ a,b\in A.
\end{equation}
Then $(A, \vartriangleright^{t}_{\l}, \vartriangleleft^{t}_{\l})$ forms an L-dendriform conformal algebra. 
Moreover, the
horizontal left-symmetric conformal algebra associated to $(A, \vartriangleright^{t}_{\l}, \vartriangleleft^{t}_{\l})$ coincides with the vertical 
left-symmetric conformal algebra $(A,\circ_{\l})$ of the original algebra 
$(A, \vartriangleright_{\l}, \vartriangleleft_{\l})$,  and the vertical left-symmetric conformal algebra associated to $(A, \vartriangleright^{t}_{\l}, \vartriangleleft^{t}_{\l})$ coincides with the horizontal left-symmetric conformal
algebra $(A,\bullet_{\l})$ of $(A, \vartriangleright_{\l}, \vartriangleleft_{\l})$. That is,
\begin{equation*}
   \bullet^{t}_{\l}=\circ_{\l} \quad\text{and}\quad \circ^{t}_{\l}=\bullet_{\l}.
\end{equation*}
\end{pro}
\begin{proof}
    It is straightforward.
\end{proof}
\begin{defi}
   Let $(A, \vartriangleright_{\l}, \vartriangleleft_{\l})$ be an L-dendriform conformal algebra. The L-dendriform conformal algebra $(A, \vartriangleright^{t}_{\l}, \vartriangleleft^{t}_{\l})$  given by Eq. \eqref{tranoperations} is called
the transpose  of $(A, \vartriangleright_{\l}, \vartriangleleft_{\l})$.
\end{defi}
\subsection{L-dendriform conformal algebras and $\mathcal{O}$-operators of left-symmetric conformal algebras}

\begin{thm}\label{thmdenonVver}
    Let $(A,\circ_{\l})$ be a left-symmetric conformal  algebra and $(V,l_A,r_A)$ be a bimodule of $A$. Suppose  $T:V\rightarrow A$ is an $\mathcal{O}$-operator associated with $(V,l_A,r_A)$.  Then 
 \begin{enumerate}
   \item [(1)] The vector space 
$V$ carries an L-dendriform conformal algebra structure defined by the operations:
\begin{equation}\label{eqLdenonV}
    u\vartriangleright^{T}_{\l}v=l_A(T(u))_{\l}(v),\quad \quad u\vartriangleleft^{T}_{\l}v=-r_A(T(u))_{\l}(v),\quad \forall~ u,v\in V.
\end{equation}
   \item [(2)] The associated vertical left-symmetric conformal algebra $(V, \circ_\l^T)$, defined by
     \begin{equation}
        u \circ_{\l}v=u \vartriangleright_{\l}^Tv-v \vartriangleleft_{-\partial-\l}^T u,\quad\forall~ u,v\in V,
      \end{equation}
   satisfies that $T:(V, \circ_\l^T)\rightarrow (A,\circ_{\l})$ is a $\mathbb{C} [\partial]$-module homomorphism of left-symmetric conformal algebras.
 
    \item [(3)] The image $T(V)=\{T(v)\;|\;v\in V\}\subset A$ is a left-symmetric conformal subalgebra of $(A,\circ_{\l})$, and inherits an L-dendriform conformal algebra structure via:
\begin{equation}\label{eqLdenonA}
    T(u)\vartriangleright_{\l}T(v)=T( u\vartriangleright^{T}_{\l}v),\quad \quad T(u)\vartriangleleft_{\l}T(v)=T(u\vartriangleleft^{T}_{\l}v),\quad \forall~ u,v\in V.
\end{equation}
   \item [(4)]  The associated vertical left-symmetric conformal algebra structure on $T(V)$ ) coincides with its left-symmetric conformal subalgebra structure induced from $(A,\circ_{\l})$, and $T$ is a $\mathbb{C} [\partial]$-module
homomorphism of L-dendriform conformal algebras from $(V, \vartriangleright_{\l}^T, \vartriangleleft_{\l}^T)$ to $(T(V), \vartriangleright_{\l}, \vartriangleleft_{\l})$.
 \end{enumerate}   
\end{thm}
\begin{proof}
For all $u, v, w\in V$, we have 
\begin{align*}
    &u\vartriangleright^{T}_{\lambda}(v\vartriangleright^{T}_{\mu}w)-(u\vartriangleright^{T}_{\l}v)\vartriangleright^{T}_{\l+\mu}w-(u\vartriangleleft^{T}_{\l}v)\vartriangleright^{T}_{\l+\mu}w-v\vartriangleright^{T}_{\mu}(u\vartriangleright^{T}_{\l}w)\\
&+(v \vartriangleright^{T}_{\mu}u)\vartriangleright^{T}_{\lambda+\mu}w+(v \vartriangleleft^{T}_{\mu}u)\vartriangleright^{T}_{\lambda+\mu}w\\
=&l_A(T(u))_{\l}l_A(T(v))_{\mu}(w)-l_A(T(v))_{\mu}l_A(T(u))_{\l}(w)-l_A(T(l_A(T(u))_{\l}(v)))_{\l+\mu}(w)\\
&-l_A(T(r_A(T(v))_{\mu}(u)))_{\l+\mu}(w)+l_A(T(l_A(T(v))_{\mu}(u)))_{\l+\mu}(w)+l_A(T(r_A(T(u))_{\l}(v)))_{\l+\mu}(w)\\
=&l_A(T(u)\circ_{\l}T(v))_{\l+\mu}(w)-l_A(T(v)\circ_{\mu}T(u))_{\l+\mu}(w)-l_A(T(l_A(T(u))_{\l}(v)))_{\l+\mu}(w)\\
&-l_A(T(r_A(T(v))_{\mu}(u)))_{\l+\mu}(w)+l_A(T(l_A(T(v))_{\mu}(u)))_{\l+\mu}(w)+l_A(T(r_A(T(u))_{\l}(v)))_{\l+\mu}(w)\\
=&l_A(T(u)\circ_{\l}T(v)-T(l_A(T(u))_{\l}(v)+r_A(T(v))_{-\partial-\l}(u)))_{\l+\mu}(w)\\
&-l_A(T(v)\circ_{\mu}T(u)-T(l_A(T(v))_{\mu}(v)+r_A(T(u))_{-\partial-\mu}(u)))_{\l+\mu}(w)=0.
\end{align*}
Similarly, we have
\begin{equation*}
   u\vartriangleright^{T}_{\lambda}(v\vartriangleleft^{T}_{\mu}w)-(u\vartriangleright^{T}_{\l}v)\vartriangleleft^{T}_{\l+\mu}w-v\vartriangleleft^{T}_{\mu}(u\vartriangleright^{T}_{\l}w)-v\vartriangleleft^{T}_{\mu}(u\vartriangleleft^{T}_{\l}w)+(v \vartriangleleft^{T}_{\mu}u)\vartriangleleft^{T}_{\lambda+\mu}w=0.
\end{equation*}
Thus $(V,\vartriangleright^{T}_{\l},\vartriangleleft^{T}_{\l})$  is an L-dendriform conformal algebra. the remaining assertions follow by direct verification. This completes the proof.
\end{proof}
\begin{cor}\label{cor3.6}
    Let $(A,\circ_{\l})$ be a left-symmetric conformal  algebra and $\mathcal{R}$ be a Rota-Baxter operator of weight zero on $A$. Then there is an L-dendriform conformal algebra structure on 
$A$ defined by the binary $\l$-operations:
\begin{equation}
    a\vartriangleright^{\mathcal{R}}_{\l}b=\mathcal{R}(a)\circ_{\l}b,\quad \quad a\vartriangleleft^{\mathcal{R}}_{\l}b=-b\circ_{-\partial-\l}\mathcal{R}(a),\quad \forall~ a,b\in A.
\end{equation}
\end{cor}
\begin{proof}
 It follows immediately from Theorem \ref{thmdenonVver} by taking $V=A$, $l_A=L_A$ and $r_A= R_A$.  
\end{proof} 

Next, we provide an example of a nontrivial Rota–Baxter operator of weight zero on a left-symmetric conformal algebra, and use it to construct an explicit example of an L-dendriform conformal algebra via Corollary \ref{cor3.6}.

\begin{ex}
Let $A=\mathbb{C}[\partial]L\oplus\mathbb{C}[\partial]W$ be a free left-symmetric conformal algebra of rank 2 with the $\lambda$-products defined by 
$$L_\lambda L=g(-\lambda)\lambda W,~L_\lambda W=W_\lambda L=W_\lambda W=0, $$
where $g(\lambda)\in \mathbb{C}[\lambda],$ $g(\lambda)\neq 0.$ This algebra was introduced in \cite[Lemma 6.13]{HongConformal2020}.

Through direct computation, we classify the nontrivial Rota–Baxter operators of weight zero on 
$A$ into the following three types:
\begin{enumerate}
  \item [{Type 1.}] $\mathcal{R} (L)=h(\partial)W,~ \mathcal{R}(W)=0$,  where $h(\partial)\in \mathbb{C}[\partial]$, $h(\partial)\neq 0$; 
  \item [{Type 2.}] $\mathcal{R}(L)=h(\partial)W,~ \mathcal{R}(W)= a W$,  where $h(\partial)\in \mathbb{C}[\partial]$, $a\in \mathbb{C}$, $a\neq 0$; 
  \item [{Type 3.}] $\mathcal{R}(L)=2a L+h(\partial)W,~ \mathcal{R}(W)= a W$,  where $h(\partial)\in \mathbb{C}[\partial]$, $a\in \mathbb{C}$, $a\neq 0$. 
\end{enumerate}

By Corollary \ref{cor3.6}, each of these Rota–Baxter operators yields an L-dendriform conformal algebra structure on $A$. However, further direct computation reveals that the operators of Type 1 and Type 2 induce only trivial structures, while the operator of Type 3 gives rise to the following nontrivial L-dendriform conformal algebra structure on $A$:
\begin{align*} 
& L\rhd_\lambda^{\mathcal{R}}L =2ag(-\lambda)\lambda W,  ~L\rhd_\lambda^{\mathcal{R}}W =  W\rhd_\lambda^{\mathcal{R}}L = W\rhd_\lambda^{\mathcal{R}}W = 0; \\
& L\lhd_\lambda^{\mathcal{R}}L =2a(\partial+\lambda)g(\partial+\lambda)W, ~ L\lhd_\lambda^{\mathcal{R}}W =  W\lhd_\lambda^{\mathcal{R}}L = W\lhd_\lambda^{\mathcal{R}}W = 0, ~~ a\in \mathbb{C}, ~ a\neq 0.
\end{align*}     
\end{ex}

Now we give the second study involving the associated horizontal left-symmetric conformal algebras.

\begin{thm}\label{thmdenonVhor}
    Let $(A,\circ_{\l})$ be a left-symmetric conformal  algebra and $(V,l_A,r_A)$ be a bimodule of $A$. Suppose  $T:V\rightarrow A$ is an $\mathcal{O}$-operator associated with $(V,l_A,r_A)$.  Then 
 \begin{enumerate}
   \item [(1)] The vector space 
$V$ carries an L-dendriform conformal algebra structure, whose $\l$-operations are defined by:
\begin{equation}
    u\vartriangleright^{T}_{\l}v=l_A(T(u))_{\l}(v),\quad \quad u\vartriangleleft^{T}_{\l}v=r_A(T(u))_{-\l-\partial}(v),\quad \forall~ u,v\in V.
\end{equation}
   \item [(2)] The associated horizontal left-symmetric conformal algebra $(V, \bullet_{\l}^T)$, defined by
     \begin{equation}
        u \bullet_{\l}v=u \vartriangleright_{\l}^Tv+u \vartriangleleft_{\l}^T v,\quad\forall~ u,v\in V,
      \end{equation}
   satisfies that $T:(V, \bullet_{\l}^T)\rightarrow (A,\circ_{\l})$ is an $\mathbb{C} [\partial]$-module homomorphism of left-symmetric conformal algebras.
 
    \item [(3)] The image $T(V)=\{T(v)\;|\;v\in V\}\subset A$ is a left-symmetric conformal subalgebra of $(A,\circ_{\l})$, and inherits an L-dendriform conformal algebra structure via:
\begin{equation}
    T(u)\vartriangleright_{\l}T(v)=T( u\vartriangleright^{T}_{\l}v),\quad \quad T(u)\vartriangleleft_{\l}T(v)=T(u\vartriangleleft^{T}_{\l}v),\quad \forall~ u,v\in V.
\end{equation}
   \item [(4)]  The associated horizontal left-symmetric conformal algebra structure on $T(V)$ ) coincides with its left-symmetric conformal subalgebra structure induced from $(A,\circ_{\l})$, and $T$ is a $\mathbb{C} [\partial]$-module
homomorphism of L-dendriform conformal algebras from $(V, \vartriangleright_{\l}^T, \vartriangleleft_{\l}^T)$ to $(T(V), \vartriangleright_{\l}, \vartriangleleft_{\l})$.
 \end{enumerate}   
\end{thm}

\begin{proof}
    It follows from a similar argument to the proof of Theorem \ref{thmdenonVver}.
\end{proof}
\begin{cor}\label{cor1}
    Let $(A,\circ_{\l})$ be a left-symmetric conformal  algebra and $\mathcal{R}$ be a Rota-Baxter operator of weight zero on $A$. Then there is an L-dendriform conformal algebra structure on $A$ defined by the binary
$\l$-operations:   
\begin{equation}
    a\vartriangleright^{\mathcal{R}}_{\l}b=\mathcal{R}(a)\circ_{\l}b,\quad \quad a\vartriangleleft^{\mathcal{R}}_{\l}b=a\circ_{\l}\mathcal{R}(b),\quad \forall~ a,b\in A.
\end{equation}
\end{cor}
\begin{proof}
 It follows immediately from Theorem \ref{thmdenonVhor} by taking $V=A$, $l_A=L_A$ and $r_A= R_A$.  
\end{proof} 

Recall that a $\mathbb{C} [\partial]$-module homomorphism $T:V\to \mathfrak{g}$ is called an $\mathcal{O}$-operator of a Lie conformal algebra $\mathfrak{g}$ associated with a representation $(V,\rho)$ if it
satisfies 
\begin{equation}
    [T(u)_{\l}T(v)]=T\Big(\rho(T(u))_{\l}v-\rho(T(v))_{-\partial-\l}u\Big),\quad \forall~ u,v\in V.
\end{equation}
In particular, an $\mathcal{O}$-operator $\mathcal{R}:\mathfrak{g} \to \mathfrak{g}$ associated with the adjoint representation $(\mathfrak{g},\mathrm{ad})$ is called a Rota-Baxter operator (of weight zero) on $\mathfrak{g}$, that is, $\mathcal{R}$ is a $\mathbb{C} [\partial]$-module homomorphism satisfying
\begin{equation}
    [\mathcal{R}(a)_{\l}\mathcal{R}(b)]=\mathcal{R}\Big([\mathcal{R}(a)_{\l}b]+[a_{\l}\mathcal{R}(b)]\Big), \quad \forall a,b\in \mathfrak{g}. 
\end{equation}
In this case,  it is known \cite{HL2} that $\mathfrak{g}$ carries a left-symmetric conformal algebra structure defined by
\begin{equation}\label{aqleftfromLie}
    a\circ_{\l}b=[\mathcal{R}(a)_{\l}b],\quad \forall a,b\in \mathfrak{g}.
\end{equation}

For brevity, in what follows we restrict our attention to the study of vertical left-symmetric conformal algebras arising from this construction. The corresponding results for horizontal left-symmetric conformal algebras can be obtained via the transposes of L-dendriform conformal algebras.
\begin{pro}
    Let $\mathfrak{g}$ be a Lie conformal algebra and $\{\mathcal{R}_1, \mathcal{R}_2\}$ be a pair of commuting $\mathcal{O}$-operators of $\mathfrak{g}$ associated with the adjoint representation $(\mathfrak{g},\mathrm{ad})$. Then there
exists an L-dendriform conformal algebra structure on $\mathfrak{g}$ defined by
\begin{equation}
     a\vartriangleright_{\l}b=[\mathcal{R}_1(\mathcal{R}_2(a))_{\l}b],\quad \quad a\vartriangleleft_{\l}b=[\mathcal{R}_2(a)_{\l}\mathcal{R}_1(b)],\quad \forall~ a,b\in \mathfrak{g}.
\end{equation}
\end{pro}

\begin{proof}
There is a left-symmetric conformal algebra structure on $\mathfrak{g}$ defined by equation \eqref{aqleftfromLie} via the $\mathcal{O}$-operator $\mathcal{R}_1$ of the Lie conformal algebra $\mathfrak{g}$ associated with $(\mathfrak{g}, \mathrm{ad})$. Since $\mathcal{R}_2$ commutes with $\mathcal{R}_1$ and is also an $\mathcal{O}$-operator of $\mathfrak{g}$ associated with $(\mathfrak{g}, \mathrm{ad})$, it follows that $\mathcal{R}_2$ is a Rota–Baxter operator of weight zero on this left-symmetric conformal algebra. The conclusion then follows from Corollary~\ref{cor1}.
\end{proof}

\begin{thm}\label{thminver}
Let $(A, \circ_\lambda)$ be a left-symmetric conformal algebra. Then there exists a compatible $L$-dendriform conformal algebra structure on $A$ such that $(A, \circ_\lambda)$ is the associated vertical left-symmetric conformal algebra if and only if there exists an invertible $\mathcal{O}$-operator of $(A, \circ_\lambda)$.
\end{thm}

\begin{proof}
Suppose there exists an invertible $\mathcal{O}$-operator $T$ of $(A, \circ_\lambda)$ associated with a bimodule $(V, l_A, r_A)$. By Theorem~\ref{thmdenonVver}, there exists an L-dendriform conformal algebra structure on $V$ given by equation~\eqref{eqLdenonV}. We then define an $L$-dendriform conformal algebra structure on $A$ via equation~\eqref{eqLdenonA} such that $T$ becomes a $\mathbb{C}[\partial]$-module isomorphism of $L$-dendriform conformal algebras, i.e.,
\begin{equation*}
    a\vartriangleright^{T}_{\l}b=T(l_A(a)_{\l}(T^{-1}(b))),\quad \quad a\vartriangleleft^{T}_{\l}b=-T(r_A(a)_{\l}(T^{-1}(b))),\quad \forall~ a,b\in A.
\end{equation*}
This structure is compatible with $(A, \circ_\lambda)$, since for all $a, b \in A$,
\begin{equation*}
  a\vartriangleright^{T}_{\l}b-b\vartriangleleft^{T}_{-\partial-\l}a=T\Big(l_A(a)_{\l}T^{-1}(b)+r_A(b)_{-\partial-\l}T^{-1}(a) \Big)=T(T^{-1}(a))\circ_{\l}T(T^{-1}(b))=a\circ_{\l}b.
\end{equation*}

Conversely, suppose $(A, \vartriangleright_\lambda, \vartriangleleft_\lambda)$ is an $L$-dendriform conformal algebra and $(A, \circ_\lambda)$ is its associated vertical left-symmetric conformal algebra. Then $(A, L_{\vartriangleright_\lambda}, -L_{\vartriangleleft_\lambda})$ is a bimodule of $(A, \circ_\lambda)$, and the identity map $\mathrm{id} : A \to A$ is an $\mathcal{O}$-operator of $(A, \circ_\lambda)$ associated with the bimodule $(A, L_{\vartriangleright_\lambda}, -L_{\vartriangleleft_\lambda})$.
\end{proof}

The following conclusion reveals the relationship between L-dendriform conformal algebras and pseudo-Hessian structures (that
is, the left-symmetric conformal algebras with a nondegenerate symmetric $2$-cocycle):

\begin{cor}
Let $(A, \circ_\lambda)$ be a left-symmetric conformal algebra equipped with a nondegenerate symmetric $2$-cocycle $B_\lambda$. Then there exists a compatible $L$-dendriform conformal algebra structure on $(A, \circ_\lambda)$ given by
\begin{equation}
\begin{aligned}
B_{\lambda+\mu}(a \vartriangleright_\lambda b, c) &= -B_\mu(b, [a_\lambda c]), \
B_{\lambda+\mu}(a \vartriangleleft_\lambda b, c) &= -B_\mu(b, c \circ_{-\partial - \lambda} a), \quad \forall a, b, c \in A,
\end{aligned}
\end{equation}
such that $(A, \circ_\lambda)$ is the associated vertical left-symmetric conformal algebra.
\end{cor}

\begin{proof}
    Let $B_{\l}$ be the nondegenerate symmetric bilinear form on $A$ defined by $B_{\l}(a,b)=\{T^{-1}(a),b\}_{\l}$, where $T:A^{*c}\to A$ is a $\mathbb{C} [\partial]$-module isomomorphism. It is straightforward to show that $T$  is an invertible $\mathcal{O}$-operator of $(A,\circ_{\l})$  associated with the bimodule $(A^{*c},L^{*}_{\circ_{\l}}-R^{*}_{\circ_{\l}},-R^{*}_{\circ_{\l}})$ . By Theorem \ref{thminver}, there is a compatible L-dendriform conformal
algebra structure on $A$ defined as follows ($a,b,c\in A$):
\begin{eqnarray*}
   B_{\l+\mu}(a\vartriangleright_{\l}b,c)&=&B_{\l+\mu}(T((L^{*}_{\circ_{\l}}-R^{*}_{\circ_{\l}})(a)T^{-1}(b)),c)=\{(L^{*}_{\circ_{\l}}-R^{*}_{\circ_{\l}})(a)T^{-1}(b),b\}_{\l+\mu}\\
   &=&- \{T^{-1}(b),(L_{\circ_{\l}}-R_{\circ_{\l}})(a)(c)\}_{\mu}=- \{T^{-1}(b),[a_{\l}c]\}_{\mu}\\
   &=&-B_{\mu}(b,[a_{\l}c]),\\
   B_{\l+\mu}(a\vartriangleleft_{\l}b,c)&=&B_{\l+\mu}(T((R^{*}_{\circ_{\l}})(a)T^{-1}(b)),c)=\{(R^{*}_{\circ_{\l}})(a)T^{-1}(b),b\}_{\l+\mu}\\
   &=&- \{T^{-1}(b),(R_{\circ_{\l}})(a)(c)\}_{\mu}=- \{T^{-1}(b),c\circ_{-\partial-\l}a\}_{\mu}\\
   &=&-B_{\mu}(b,c\circ_{-\partial-\l}a).
\end{eqnarray*}
Hence $(A,\circ_{\l})$ is the associated vertical left-symmetric conformal algebra, and the result follows.
\end{proof}

The following result presents a construction of (symmetric) solutions to the conformal $\mathcal{S}$-equation in a class of left-symmetric conformal algebras arising from L-dendriform conformal algebras:

\begin{thm}
Let $(A, \vartriangleright_{\l}, \vartriangleleft_{\l})$ be a free  L-dendriform conformal algebra of finite rank,  and let  $(A,\circ_{\l})$ and $(A,\bullet_{\l})$ be the associated vertical and horizontal left-symmetric conformal
algebras, respectively. Then
\begin{equation}
r = \sum_{i=1}^{n} \left( e_i \otimes e^*_i + e^*_i \otimes e_i \right)
\tag{61}
\end{equation}
is a symmetric solution of the conformal $S$-equation in the left-symmetric conformal algebras
$$
A \ltimes_{L^{*}_{\vartriangleright_{\l}}+L^{*}_{\vartriangleleft_{\l}},L^{*}_{\vartriangleleft_{\l}}} A^{*c}~~\text{and} ~~
A \ltimes_{L^{*}_{\vartriangleright_{\l}}-R^{*}_{\vartriangleleft_{\l}},-R^{*}_{\vartriangleleft_{\l}}} A^{*c}, $$ where \( \{e_1, \dots, e_n\} \) is a \( \mathbb{C}[\partial] \)-basis of \( A \) and \( \{e^*_1, \dots, e^*_n\} \) is the dual \( \mathbb{C}[\partial] \)-basis of \( A^{*c} \).   
\end{thm}

\begin{proof}
By Proposition \ref{horver} and Corollary \ref{corhorver}, the map $\mathrm{id} : A \rightarrow A $ is an \( \mathcal{O} \)-operator  associated with the bimodules $(A,L_{\vartriangleright_{\l}},-L_{\vartriangleleft_{\l}})$, and $(A,L_{\vartriangleright_{\l}},R_{\vartriangleleft_{\l}})$, respectively. The result then follows directly from Theorem~\ref{solS-equa}.
\end{proof}
\subsection{Relationships with dendriform conformal algebras and quadri conformal algebras}
\begin{defi}  (\cite{HB})
    A dendriform conformal algebra is a triple $
(A, \succ_\lambda, \prec_\lambda)$ consisting of a $ \mathbb{C}[\partial] $-module $A$ and two $\lambda$-multiplications
\[
\succ_\lambda, ~\prec_\lambda : A \times A \to A[\lambda],
\]
which are conformal sesquilinear maps and satisfy the following identities:
\begin{align*}
a \succ_\lambda (b \succ_\mu c) &= (a \succ_\lambda b + a \prec_\lambda b) \succ_{\lambda+\mu} c,\\ 
 (a \prec_\lambda b) \prec_{\lambda+\mu} c &= a \prec_\lambda (b \succ_\mu c + b \prec_\mu c),\\ 
(a \succ_\lambda b) \prec_{\lambda+\mu} c &= a \succ_\lambda (b \prec_\mu c),
\end{align*}
for all $a, b, c \in A.$
\end{defi}
\begin{pro}\label{denLden}
    Any dendriform conformal algebra $
(A, \succ_\lambda, \prec_\lambda)$ carries a natural 
L-dendriform conformal algebra structure under the operations:  $$a\vartriangleright_{\l}b=a\succ_\lambda b,~~a\vartriangleleft_{\l}b=a\prec_\lambda b.$$
\end{pro}
\begin{proof} Let $(A, \succ_\lambda, \prec_\lambda)$ be a dendriform conformal algebra.
The identities \eqref{eq1:L-dendriform conformal1} and \eqref{eq2:L-dendriform conformal2} that define an 
L-dendriform conformal algebra are both satisfied trivially in this case, since both sides of each identity evaluate to zero under the given operations. This follows directly from the defining axioms of a dendriform conformal algebra. Hence, the claim holds.
\end{proof}

\begin{rmk}
In this sense, associative conformal algebras are precisely those left-symmetric conformal algebras that are associative (i.e., whose associator vanishes), whereas dendriform conformal algebras are precisely those L-dendriform conformal algebras for which the ``generalized associators'' (see Remark~\ref{genass}) vanish.
\end{rmk}


By Propositions \ref{denleft} and \ref{denLden}, the following result holds:

\begin{cor}\label{cor3.14}
  Let $(A, \succ_\lambda, \prec_\lambda)$ be a dendriform conformal algebra. Then 
  \begin{itemize}
\item[(1)] The binary $\l$-operation defined by Eq.~\eqref{eqbullet} yields a left-symmetric conformal algebra (in fact, an associative conformal algebra).
\item[(2)] The binary $\l$-operation defined by Eq.~\eqref{eqcirc} yields a left-symmetric conformal algebra.
\end{itemize} 
\end{cor}

\begin{defi} (\cite{AW})
     A  quadri conformal algebra is a 5-tuple $(A, \nwarrow_{\l},
\swarrow_{\l}, \nearrow_{\l}, \searrow_{\l})$
consisting of a $ \mathbb{C}[\partial] $-module $A$ and four $\lambda$-multiplications $$\nwarrow_{\l},
\swarrow_{\l}, \nearrow_{\l}, \searrow_{\l}\, : A \times A\rightarrow A[\lambda]$$ which are conformal sesquilinear maps and satisfy the following identities for all $a, b, c \in A$:
\small{\begin{align}
(a \nwarrow_{\l}b) \nwarrow_{\l+\mu} c &= a \nwarrow_{\l}(b\ast_{\mu}c),\quad  (a \nearrow_{\l}b) \nwarrow_{\l+\mu} c = a  \nearrow_{\l}(b\prec_{\mu}c),\quad (a \wedge_{\l}b) \nearrow_{\l+\mu} c = a  \nearrow_{\l}(b\succ_{\mu}c),\\ 
(a \swarrow_{\l}b) \nwarrow_{\l+\mu} c &= a \swarrow_{\l}(b\wedge_{\mu}c),\quad  (a \vee_{\l}b) \nearrow_{\l+\mu} c = a  \searrow_{\l}(b\nearrow_{\mu}c),\quad (a \searrow_{\l}b) \nwarrow_{\l+\mu} c = a  \searrow_{\l}(b\nwarrow_{\mu}c),\\ 
(a \prec_{\l}b) \swarrow_{\l+\mu} c &= a \swarrow_{\l}(b\vee_{\mu}c),\quad  (a \ast_{\l}b) \searrow_{\l+\mu} c = a  \searrow_{\l}(b\searrow_{\mu}c),\quad (a \succ_{\l}b) \swarrow_{\l+\mu} c = a  \searrow_{\l}(b\swarrow_{\mu}c),
\end{align}}
where
\begin{align}
a \vee_{\l}b &= a \swarrow_{\l}b+a\searrow_{\l} b,\quad  a \wedge_{\l}b = a \nwarrow_{\l}b+a\nearrow_{\l} b,\label{dendri1}\\ 
a \succ_{\l}b &= a \nwarrow_{\l}b+a\swarrow_{\l} b,\quad  a \prec_{\l}b = a \nearrow_{\l}b+a\searrow_{\l} b,\label{dendri2}\\ 
a \ast_{\l}b &= a \swarrow_{\l}b+a\searrow_{\l} b+ a \nwarrow_{\l}b+a\nearrow_{\l} b=a \vee_{\l}b+a \wedge_{\l}b=a \succ_{\l}b+ a \prec_{\l}b.\label{dendri3}
\end{align}
\end{defi}
\begin{pro}\label{3.15}(\cite{AW})
 Let $(A, \nwarrow_{\l},
\swarrow_{\l}, \nearrow_{\l}, \searrow_{\l})$ be a quadri conformal algebra.
\begin{enumerate}
\item[(1)] $
(A, \succ_\lambda, \prec_\lambda)$ forms a dendriform conformal algebra, with the $\lambda$-operations given in Eq. \eqref{dendri1}.
    \item[(2)] $
(A, \vee_\lambda, \wedge_\lambda)$  also forms a dendriform conformal algebra, with the $\lambda$-operations given by Eq. \eqref{dendri2}.
\item[(3)] $
(A, \ast_\lambda)$ is an associative conformal algebra,  where the $\lambda$-multiplication $\ast_\lambda$ is given by Eq. \eqref{dendri3}.
\end{enumerate}
\end{pro}

There is a further connection between quadri conformal algebras and L-dendriform conformal algebras, which can be stated as follows: 
\begin{pro}\label{3.16}
  Let $(A, \nwarrow_{\l},
\swarrow_{\l}, \nearrow_{\l}, \searrow_{\l})$ be  quadri conformal algebra. The $\lambda$-multiplications defined by 
\begin{equation}
a \vartriangleright_{\l}b=a \searrow_{\l}b-b\nwarrow_{-\partial -\l}a, \quad  a\vartriangleleft_{\l}b=  a \nearrow_{\l}b-b\swarrow_{-\partial -\l}a, 
\end{equation}
yield an L-dendriform conformal algebra structure on $A$.
\end{pro}
\begin{proof}
    The result follows by a direct verification.
\end{proof}

\begin{cor}
   Let $(A, \nwarrow_{\l},
\swarrow_{\l}, \nearrow_{\l}, \searrow_{\l})$ be a quadri conformal algebra. For all $a,b\in A$:
\begin{enumerate}
    \item[(1)] The $\lambda$-multiplication given by 
\begin{equation}
a \circ_{\l}b=a \searrow_{\l}b+a\swarrow_{\l}b-b\nwarrow_{-\partial -\l}a-  b \nearrow_{-\partial -\l}a=a \vartriangleright_{\l}b-b\vartriangleleft_{-\partial -\l}a=a \vee_{\l}b-b\wedge_{-\partial -\l}a, 
\end{equation}
defines a left-symmetric conformal algebra $(A,\circ_{\l})$.
\item[(2)] The $\lambda$-multiplication given by 
\begin{equation}
a \bullet_{\l}b=a \searrow_{\l}b+a \nearrow_{\l}b-b\nwarrow_{-\partial -\l}a -b\swarrow_{-\partial -\l}a=a \vartriangleright_{\l}b+a\vartriangleleft_{\l}b=a \succ_{\l}b-b\prec_{-\partial -\l}a, 
\end{equation}
defines a left-symmetric conformal algebra $(A,\bullet_{\l})$.
\item[(3)] The $\lambda$-multiplication given by 
\small{\begin{equation}
[a_{\l}b]=a \swarrow_{\l}b+a\searrow_{\l} b+ a \nwarrow_{\l}b+a\nearrow_{\l} b-(b \swarrow_{-\partial -\l}a+b\searrow_{-\partial -\l} a+ b \nwarrow_{-\partial -\l}a+b\nearrow_{-\partial -\l} a, 
\end{equation}}
defines a Lie conformal algebra $(\mathfrak{g}(A),[\c_{\l}\c])$.
\end{enumerate}
\end{cor}
\begin{proof}
Claims (1) and (2) follow from Propositions \ref{3.15} and \ref{3.16} along with Corollary \ref{cor3.14}. Claim (3) follows from Propositions \ref{3.15} and \ref{LEFTLIE}.
\end{proof}

Summarizing the above study in this subsection, we have the following commutative diagram:\\

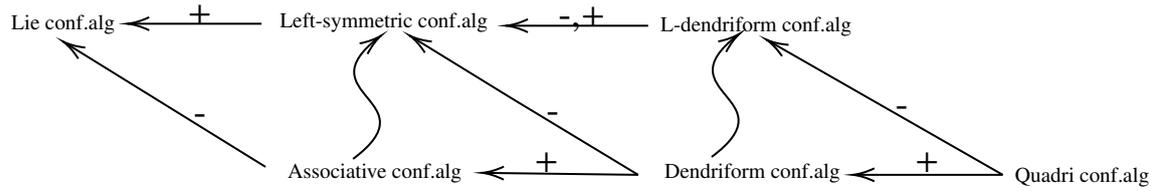
\tikzset{every picture/.style={line width=0.75pt}} 
\begin{tikzpicture}[x=0.75pt,y=0.75pt,yscale=-1,xscale=1]

\draw    (327.33,91) -- (258.33,91) ;
\draw [shift={(256.33,91)}, rotate = 360] [color={rgb, 255:red, 0; green, 0; blue, 0 }  ][line width=0.75]    (10.93,-3.29) .. controls (6.95,-1.4) and (3.31,-0.3) .. (0,0) .. controls (3.31,0.3) and (6.95,1.4) .. (10.93,3.29)   ;
\draw    (134.33,90) -- (69.33,90) ;
\draw [shift={(67.33,90)}, rotate = 360] [color={rgb, 255:red, 0; green, 0; blue, 0 }  ][line width=0.75]    (10.93,-3.29) .. controls (6.95,-1.4) and (3.31,-0.3) .. (0,0) .. controls (3.31,0.3) and (6.95,1.4) .. (10.93,3.29)   ;
\draw    (504.33,166) -- (431.33,166) ;
\draw [shift={(429.33,166)}, rotate = 360] [color={rgb, 255:red, 0; green, 0; blue, 0 }  ][line width=0.75]    (10.93,-3.29) .. controls (6.95,-1.4) and (3.31,-0.3) .. (0,0) .. controls (3.31,0.3) and (6.95,1.4) .. (10.93,3.29)   ;
\draw    (322.33,166) -- (246.33,165.03) ;
\draw [shift={(244.33,165)}, rotate = 0.73] [color={rgb, 255:red, 0; green, 0; blue, 0 }  ][line width=0.75]    (10.93,-3.29) .. controls (6.95,-1.4) and (3.31,-0.3) .. (0,0) .. controls (3.31,0.3) and (6.95,1.4) .. (10.93,3.29)   ;
\draw    (136.33,162) -- (37.03,100.06) ;
\draw [shift={(35.33,99)}, rotate = 31.95] [color={rgb, 255:red, 0; green, 0; blue, 0 }  ][line width=0.75]    (10.93,-3.29) .. controls (6.95,-1.4) and (3.31,-0.3) .. (0,0) .. controls (3.31,0.3) and (6.95,1.4) .. (10.93,3.29)   ;
\draw    (504.33,166) -- (388.07,99) ;
\draw [shift={(386.33,98)}, rotate = 29.95] [color={rgb, 255:red, 0; green, 0; blue, 0 }  ][line width=0.75]    (10.93,-3.29) .. controls (6.95,-1.4) and (3.31,-0.3) .. (0,0) .. controls (3.31,0.3) and (6.95,1.4) .. (10.93,3.29)   ;
\draw    (322.33,166) -- (258.33,127) -- (210.03,97.05) ;
\draw [shift={(208.33,96)}, rotate = 31.8] [color={rgb, 255:red, 0; green, 0; blue, 0 }  ][line width=0.75]    (10.93,-3.29) .. controls (6.95,-1.4) and (3.31,-0.3) .. (0,0) .. controls (3.31,0.3) and (6.95,1.4) .. (10.93,3.29)   ;
\draw    (180,158) .. controls (219.6,128.3) and (155.31,126.04) .. (192.83,96.89) ;
\draw [shift={(194,96)}, rotate = 143.13] [color={rgb, 255:red, 0; green, 0; blue, 0 }  ][line width=0.75]    (10.93,-3.29) .. controls (6.95,-1.4) and (3.31,-0.3) .. (0,0) .. controls (3.31,0.3) and (6.95,1.4) .. (10.93,3.29)   ;
\draw    (359,157) .. controls (398.6,127.3) and (334.31,125.04) .. (371.83,95.89) ;
\draw [shift={(373,95)}, rotate = 143.13] [color={rgb, 255:red, 0; green, 0; blue, 0 }  ][line width=0.75]    (10.93,-3.29) .. controls (6.95,-1.4) and (3.31,-0.3) .. (0,0) .. controls (3.31,0.3) and (6.95,1.4) .. (10.93,3.29)   ;

\draw (8,84) node [anchor=north west][inner sep=0.75pt]  [font=\scriptsize] [align=left] {Lie conf.alg};
\draw (142,83) node [anchor=north west][inner sep=0.75pt]  [font=\scriptsize] [align=left] {Left-symmetric conf.alg};
\draw (332,85) node [anchor=north west][inner sep=0.75pt]  [font=\scriptsize] [align=left] {L-dendriform conf.alg};
\draw (509,161) node [anchor=north west][inner sep=0.75pt]  [font=\scriptsize] [align=left] {Quadri conf.alg};
\draw (334,159) node [anchor=north west][inner sep=0.75pt]  [font=\scriptsize] [align=left] {Dendriform conf.alg};
\draw (146,159) node [anchor=north west][inner sep=0.75pt]  [font=\scriptsize] [align=left] {Associative conf.alg};
\draw (96,79) node [anchor=north west][inner sep=0.75pt]  [font=\Large] [align=left] {+};
\draw (99,134) node [anchor=north west][inner sep=0.75pt]  [font=\Large] [align=left] {\mbox{-}};
\draw (281,80) node [anchor=north west][inner sep=0.75pt]  [font=\Large] [align=left] {\mbox{-},+};
\draw (449.33,130) node [anchor=north west][inner sep=0.75pt]  [font=\Large] [align=left] {\mbox{-}};
\draw (269,153) node [anchor=north west][inner sep=0.75pt]  [font=\Large] [align=left] {+};
\draw (459,153) node [anchor=north west][inner sep=0.75pt]  [font=\Large] [align=left] {+};
\draw (275,133) node [anchor=north west][inner sep=0.75pt]  [font=\Large] [align=left] {\mbox{-}};

\end{tikzpicture}
\\
where ``$\leadsto $'' means inclusion, ``$+$'' means the binary $\l$-operation $x \circ^{1}_{\l} y + x \circ^{2}_{\l} y$ and ``$-$'' means the binary $\l$-operation $x \circ^{1}_{\l} y - y \circ^{2}_{-\partial-\l} x$.

\section{Compatible $\mathcal{O}$-operators and compatible L-dendriform conformal 
algebras }\label{s4}

In this section, we explore the connections between Nijenhuis operators and certain compatible structures such as compatible $\mathcal{O}$-operators and compatible L-dendriform conformal algebras through the framework of $\mathcal{O}$-operators on left-symmetric conformal algebras.

Recall that a Nijenhuis operator $N$ on a left-symmetric conformal algebra $(A, \circ_{\l})$ is a $\mathbb{C}[\partial]$-module homomorphism $N : A \to A$ satisfying
\begin{equation}
N(a) \circ_{\l} N(b) = N\big(N(a) \circ_{\l} b + a \circ_{\l} N(b) - N(a \circ_{\l} b)\big), \quad \forall~ a, b \in A.
\end{equation}
For instance, the identity map $\mathrm{id}$ is a Nijenhuis operator on $A$.

Let $N$ be a Nijenhuis operator on a left-symmetric conformal algebra $(A, \circ_{\l})$. We define a new $\l$-product 
\begin{equation}
a\circ^{N}_{\l}b=N(a)\circ_{\l}b+a\circ_{\l}N(b)-N(a\circ_{\l}b),\quad \forall a,b\in A.
\end{equation} 
Then $(A,\circ^{N}_{\l})$ also carries the structure of a left-symmetric conformal algebra, and $N$ is a homomorphism from  $(A,\circ^{N}_{\l})$ 
to $(A,\circ_{\l})$.

\begin{defi}
Let $A$ be a left-symmetric conformal algebra and $\mathcal{R} : A \to A$ a $\mathbb{C}[\partial]$-module homomorphism. For $q \in \mathbb{C}$, if $\mathcal{R}$ satisfies
\begin{equation}
\mathcal{R}(a) \circ_{\l} \mathcal{R}(b) = \mathcal{R}\Big( \mathcal{R}(a) \circ_{\l} b + a \circ_{\l} \mathcal{R}(b) + q a \circ_{\l} b \Big), \quad \forall a, b \in A,
\end{equation}
then $\mathcal{R}$ is called a \textbf{Rota–Baxter operator of weight $q$} on $A$.
\end{defi}

The following proposition establishes a close relationship between Nijenhuis operators and Rota–Baxter operators; its proof follows by direct verification.

\begin{pro}\label{nijenhuispro}
Let $N: A \rightarrow A$ be a $\mathbb{C}[\partial]$-module homomorphism over a  left-symmetric conformal algebra $A$.

\begin{itemize}
    \item[(i)] If $N^2 = 0$, then $N$ is a Nijenhuis operator if and only if $N$ is a Rota-Baxter operator of weight $0$.
    \item[(ii)] If $N^2 = N$, then $N$ is a Nijenhuis operator if and only if $N$ is a Rota-Baxter operator of weight $-1$.
    \item[(iii)] If $N^2 = \mathrm{id}$, then $N$ is a Nijenhuis operator if and only if $N \pm \mathrm{id}$ is a Rota-Baxter operator of weight $\mp 2$.
\end{itemize}
    
\end{pro}
\begin{ex}
 Let $T: V \rightarrow A$ be an $\mathcal{O}$-operator on a left-symmetric conformal algebra $A$ with respect to a bimodule $V$ of $A$. By Proposition \ref{caraco-operator}, the lift $\hat{T}$ is a Rota-Baxter operator of weight $0$ on $A \oplus_0 V$. Clearly, $\hat{T}^2 = 0$. Hence, by Proposition \ref{nijenhuispro} (i), $\hat{T}$ is a Nijenhuis operator on $A \oplus_0 V$.  
\end{ex}
\begin{defi}
Let $(A,\circ_{\l})$ be a left-symmetric conformal algebra and $(V;l_A,r_A)$ be a bimodule of $A$.  Let $T_1, T_2: V \rightarrow A$ be two $\mathcal{O}$-operators.  
If for all $k_1, k_2 \in \mathbb{C}$, the linear combination $k_1 T_1 + k_2 T_2$ is still an $\mathcal{O}$-operator, then $T_1$ and $T_2$ are called {\bf compatible}.   
\end{defi}

\begin{pro}
Let $T_1, T_2: V \rightarrow A$ be two $\mathcal{O}$-operators on a left-symmetric conformal algebra $(A,\circ_{\l})$ associated to the representation 
$(V;l_A,r_A)$. Then $T_1$ and $T_2$ are compatible if and only if the following equation holds for all $u,v \in V:$
\begin{align}\label{compaoper1}
T_1(u)\circ_{\l} T_2(v) + T_2(u) \circ_{\l} T_1(v) =\; & T_1\big(l_A(T_2(u))_\lambda v+r_A(T_2(v))_{-\partial-\lambda}u\big) \notag \\
& + T_2\big(l_A(T_1(u))_\lambda v+r_A(T_1(v))_{-\partial-\lambda}u\big).
\end{align}
    
\end{pro} 
\begin{proof}
 This follows by a direct calculation using the definitions of 
$\mathcal{O}$-operators and the compatibility condition.   
\end{proof}

Using an $\mathcal{O}$-operator and a Nijenhuis operator, one can construct a pair of compatible $\mathcal{O}$-operators.

\begin{pro}\label{compatible1}
Let $T: V \rightarrow A$ be an $\mathcal{O}$-operator on a left-symmetric conformal algebra  $(A,\circ_{\l})$ with respect to the representation $(V,l_A,r_A)$ and let $N$ be a Nijenhuis operator on $A$. Then $N T$ is also an $\mathcal{O}$-operator on $(A,\circ_{\l})$ if and only if for all $u, v \in V$, the following equation holds:
\begin{align}\label{compaoper2}
N\big(NT(u)\circ_{\l} T(v) + T(u) \circ_{\l} NT(v)\big) =\; & N\Big(T(l_A(NT(u))_\lambda v+r_A(NT(v))_{-\partial-\lambda}u) \notag \\
& + NT(l_A(T(u))_\lambda v+r_A(T(v))_{-\partial-\lambda}u)\Big). 
\end{align}
Furthermore, if $N$ is invertible, then $T$ and $N T$ are compatible. 
\end{pro} 
\begin{proof}
For all $u, v \in V$, since $N$ is a Nijenhuis operator, we have
\[
NT(u)\circ_{\l} NT(v) = N\big(NT(u) \circ_{\l} T(v) + T(u)\circ_{\l} NT(v)\big) - N^2(T(u) \circ_{\l} T(v)).
\]
Note that since $T$ is an $\mathcal{O}$-operator, 
\[
T(u) \circ_{\l} T(v) = T\big(l_A(T(u))_\lambda v+r_A(T(v))_{-\partial-\lambda}u\big).
\]
Therefore,
\[
NT(u) \circ_{\l} NT(v) = NT\big(l_A(NT(u))_\lambda v+r_A(NT(v))_{-\partial-\lambda}u\big)
\]
if and only if the identity in Eq. \eqref{compaoper2} is satisfied.

Now, if $NT$ is an $\mathcal{O}$-operator and $N$ is invertible, it follows that
\[
NT(u) \circ_{\l} T(v) + T(u) \circ_{\l} NT(v) = T\big(l_A(NT(u))_\lambda v+r_A(NT(v))_{-\partial-\lambda}u\big) + NT\big(l_A(T(u))_\lambda v+r_A(T(v))_{-\partial-\lambda}u\big),
\]
which is exactly the compatibility condition between $T$ and $N T$.  
\end{proof}
A pair of compatible $\mathcal{O}$-operators can also give rise to a Nijenhuis operator under certain conditions.
\begin{pro}\label{compatible2}
Let $T_1, T_2: V \rightarrow A$ be two $\mathcal{O}$-operators on a left-symmetric conformal algebra $(A,\circ_{\l})$ with respect to the representation $(V,l_A,r_A)$ and suppose $T_2$ is invertible. If $T_1$ and $T_2$ are compatible, then $N = T_1  T_2^{-1}$ is a Nijenhuis operator on $A$.  
\end{pro} 
\begin{proof}
For any $a,b \in A$, there exist $u, v \in V$ such that $a= T_2(u) $, $b=T_2(v) $. Then $N = T_1 T_2^{-1}$ is a Nijenhuis operator if and only if 
\begin{align*}
N(a) \circ_{\l} N(b) = N\big(N(a) \circ_{\l} b + a \circ_{\l} N(b)\big) - N^2(a \circ_{\l} b).
\end{align*}
That is, 
\[
NT_2(u) \circ_{\l} NT_2(v) = N\big(NT_2(u) \circ_{\l} T_2(v) + T_2(u)\circ_{\l} NT_2(v)\big) - N^2\big(T_2(u) \circ_{\l} T_2(v)\big).
\]
Since $T_1 = N  T_2$ is an $\mathcal{O}$-operator, the left-hand side becomes
\[
NT_2\big(l_A(NT_2(u))_\lambda v+r_A(NT_2(v))_{-\partial-\lambda}u\big).
\]
On the other hand,
since $T_2$ is an $\mathcal{O}$-operator and is compatible with $T_1 = N  T_2$, we have
\begin{align*}
NT_2(u) \circ_{\l} T_2(v) + T_2(u) \circ_{\l} NT_2(v)
&= T_2\big(l_A(NT_2(u))_\lambda v+r_A(NT_2(v))_{-\partial-\lambda}u\big) \\
&\quad + NT_2\big(l_A(T_2(u))_\lambda v+r_A(T_2(v))_{-\partial-\lambda}u\big) \\
&= T_2\big(l_A(NT_2(u))_\lambda v+r_A(NT_2(v))_{-\partial-\lambda}u\big) + N\big(T_2(u) \circ_{\l} T_2(v)\big).
\end{align*}
Applying $N$ to both sides of this identity yields the required Nijenhuis condition.   
\end{proof}

By Propositions \ref{compatible1} and \ref{compatible2}, we obtain the following corollary:

\begin{cor}
 Let $T_1, T_2: V \rightarrow A$ be two $\mathcal{O}$-operators on a left-symmetric conformal algebra $(A,\circ_{\l})$ with respect to the representation $(V,l_A,r_A)$ and suppose both $T_1$ and $T_2$ are invertible. Then $T_1$ and $T_2$ are compatible if and only if $N = T_1 T_2^{-1}$ is a Nijenhuis operator.    
\end{cor}
Since a Rota--Baxter operator of weight zero is an $\mathcal{O}$-operator on a left-symmetric conformal algebra associated to the adjoint bimodule, we further deduce:

\begin{cor}
 Let $(A,\circ_{\l})$ be a left-symmetric conformal algebra, and suppose both $\mathcal{R}_1$ and $\mathcal{R}_2$ are two invertible Rota-Baxter operators of weight zero on $A$. Then $\mathcal{R}_1$ and $\mathcal{R}_2$ are compatible in the sense that every linear combination of $\mathcal{R}_1$ and $\mathcal{R}_2$ remains a Rota-Baxter operator of weight zero if and only if $N = \mathcal{R}_1 \mathcal{R}_2^{-1}$ is a Nijenhuis operator.   
\end{cor} 

\begin{rmk}
  Note that the inverse of an invertible Rota-Baxter operator is a derivation on a left-symmetric conformal algebra. However, it is interesting to see that a similar role of Nijenhuis operators as above is not available for derivations, since any linear combination of two derivations is always a derivation due to the fact that the set of all derivations is a linear space.  
\end{rmk} 

The notion of an L-dendriform conformal algebra was introduced in this paper as the natural algebraic structure underlying an $\mathcal{O}$-operator on a left-symmetric conformal algebra. To conclude this section, we define the concept of compatible L-dendriform conformal algebras and demonstrate that compatible $\mathcal{O}$-operators on a left-symmetric conformal algebra naturally induce such structures.

\begin{defi}
 Two L-dendriform conformal algebras $(A, \vartriangleright^{1}_{\l}, \vartriangleleft^{1}_{\l})$ and $(A, \vartriangleright^{2}_{\l}, \vartriangleleft^{2}_{\l})$ are said to be compatible if for all $k_1, k_2 \in \mathbb{C}$, the structure $(A, k_1 \vartriangleright^{1}_{\l} + k_2 \vartriangleright^{2}_{\l}, k_1 \vartriangleleft^{1}_{\l} + k_2 \vartriangleleft^{2}_{\l})$ is again an L-dendriform conformal algebra. We denote it by
$
(A, \vartriangleright^{1}_{\l}, \vartriangleleft^{1}_{\l}; \vartriangleright^{2}_{\l}, \vartriangleleft^{2}_{\l}).$  
\end{defi} 

\begin{lem}
 Two L-dendriform conformal algebras $(A, \vartriangleright^{1}_{\l}, \vartriangleleft^{1}_{\l})$ and $(A, \vartriangleright^{2}_{\l}, \vartriangleleft^{2}_{\l})$ are compatible if and only if for all $a, b, c \in A$, the following equations hold:
\begin{align*}
a \vartriangleright^{1}_{\l} (b \vartriangleright^{2}_{\mu} c) + a \vartriangleright^{2}_{\l} (b \vartriangleright^{1}_{\mu} c)
&= (a \vartriangleright^{1}_{\l} b) \vartriangleright^{2}_{\l+\mu} c + (a \vartriangleright^{2}_{\l} b) \vartriangleright^{1}_{\l+\mu} c + (a \vartriangleleft^{1}_{\l} b) \vartriangleright^{1}_{\l+\mu} c + (a \vartriangleleft^{2}_{\l}  b) \vartriangleright^{1}_{\l+\mu} c \\
&\quad + b \vartriangleright^{1}_{\mu} (a \vartriangleright^{2}_{\l}c) + b \vartriangleright^{2}_{\mu} (a \vartriangleright^{1}_{\l} c)  - (b \vartriangleleft^{1}_{\mu} a) \vartriangleright^{2}_{\l+\mu} c - (b \vartriangleleft^{2}_{\mu} a) \vartriangleright^{1}_{\l+\mu} c \\
&\quad - (b \vartriangleright^{1}_{\mu} a) \vartriangleright^{2}_{\l+\mu} c - (b \vartriangleright^{2}_{\mu}a) \vartriangleright^{1}_{\l+\mu} c,\\
a \vartriangleright^{1}_{\l} (b \vartriangleleft^{2}_{\mu} c) + a \vartriangleright^{2}_{\l} (b \vartriangleleft^{1}_{\mu} c)
&= (a \vartriangleright^{1}_{\l} b) \vartriangleleft^{2}_{\l+\mu} c + (a \vartriangleright^{2}_{\l} b) \vartriangleleft^{1}_{\l+\mu} c + b \vartriangleleft^{1}_{\mu} (a \vartriangleright^{2}_{\l} c) + b \vartriangleleft^{2}_{\mu} (a \vartriangleright^{1}_{\l} c) \\
&\quad + b \vartriangleleft^{2}_{\mu} (a \vartriangleleft^{1}_{\l} c) + b \vartriangleleft^{1}_{\mu} (a \vartriangleleft^{2}_{\l} c) - (b \vartriangleleft^{1}_{\mu} a)\vartriangleleft^{2}_{\l+\mu} c - (b \vartriangleleft^{2}_{\mu} a) \vartriangleleft^{1}_{\l+\mu} c.
\end{align*}  
\end{lem} 
\begin{proof}
    It is straightforward.
\end{proof}
\begin{pro} Let $(A,\circ_{\l})$ be a left-symmetric conformal algebra and $(V,l_A,r_A)$ a bimodule of $A$.
 Suppose $T_1, T_2 : V \rightarrow A$ are two compatible $\mathcal{O}$-operators.  Then the triples $(V, \vartriangleright^{1}_{\l}, \vartriangleleft^{1}_{\l})$ and $(V, \vartriangleright^{2}_{\l}, \vartriangleleft^{2}_{\l})$ form compatible L-dendriform conformal algebras, where
   \begin{eqnarray*}
      &&u \vartriangleright^{1}_{\l} v = l_A(T_1(u))_{\l}v, \quad u \vartriangleleft^{1}_{\l} v = -r_A(T_1(u))_{\l}v, \\
      &&u \vartriangleright^{2}_{\l} v = l_A(T_2(u))_{\l}v, \quad u \vartriangleleft^{2}_{\l} v = -r_A(T_2(u))_{\l}v, \quad \forall~ u, v \in V.
   \end{eqnarray*}

Moreover, if both $T_1$ and $T_2$ are invertible, then $(A, \vartriangleright^{1}_{\l}, \vartriangleleft^{1}_{\l})$ and $(A, \vartriangleright^{2}_{\l}, \vartriangleleft^{2}_{\l})$ are also compatible L-dendriform conformal algebras, where
\begin{eqnarray*}
 &&a \vartriangleright^{1}_{\l} b = T_1 l_A(a)_{\l} T_1^{-1}(b), \quad a \vartriangleleft^{1}_{\l} b = -T_1 r_A(a)_{\l} T_1^{-1}(b),\\
&&a \vartriangleright^{2}_{\l} b = T_2 l_A(a)_{\l} T_2^{-1}(b), \quad a \vartriangleleft^{2}_{\l} b = -T_2 r_A(a)_{\l} T_2^{-1}(b), \quad \forall~ a, b \in A.   
\end{eqnarray*}
\end{pro}

\begin{proof}
 This follows immediately from Theorems \ref{thmdenonVhor} and \ref{thminver}, together with the definitions of compatible $\mathcal{O}$-operators and compatible L-dendriform conformal algebras.   
\end{proof} 

\begin{cor}
Let $T : V \rightarrow A$ be two $\mathcal{O}$-operators on a left-symmetric conformal algebra $(A,\circ_{\l})$ with respect to the representation $(V,l_A,r_A)$. Suppose there exists an invertible Nijenhuis operator $N$ such that $N T$ is also an $\mathcal{O}$-operator on $A$. Then the triples $(V, \vartriangleright^{1}_{\l}, \vartriangleleft^{1}_{\l})$ and $(V, \vartriangleright^{2}_{\l}, \vartriangleleft^{2}_{\l})$ form compatible L-dendriform conformal algebras, where
 \begin{eqnarray*}
      &&u \vartriangleright^{1}_{\l} v = l_A(T(u))_{\l}v, \quad u \vartriangleleft^{1}_{\l} v = -r_A(T(u))_{\l}v, \\
      &&u \vartriangleright^{2}_{\l} v = l_A(NT(u))_{\l}v, \quad u \vartriangleleft^{2}_{\l} v = -r_A(NT(u))_{\l}v, \quad \forall u, v \in V.
   \end{eqnarray*}
Furthermore, if $T$ is also invertible, then $(A, \vartriangleright^{1}_{\l}, \vartriangleleft^{1}_{\l})$ and $(A, \vartriangleright^{2}_{\l}, \vartriangleleft^{2}_{\l})$ are compatible L-dendriform conformal algebras, with the operations defined for all $ a, b \in A$ by 
\begin{eqnarray*}
 &&a \vartriangleright^{1}_{\l} b = T l_A(a)_{\l} T^{-1}(b), \quad\quad a \vartriangleleft^{1}_{\l} b = -T r_A(a)_{\l} T^{-1}(b),\\
&&a \vartriangleright^{2}_{\l} b = NT l_A(a)_{\l} (NT)^{-1}(b), \quad a \vartriangleleft^{2}_{\l} b = -NT r_A(a)_{\l} (NT)^{-1}(b).  
\end{eqnarray*}   
\end{cor}

\begin{ex}
Let $\mathcal R$ be a Rota-Baxter operator of weight zero on a left-symmetric conformal algebra $(A,\circ_{\l})$. Suppose that there exists an invertible Nijenhuis operator $N$ such that $N \mathcal R$ is also a Rota-Baxter operator of weight zero. Then $(A, \vartriangleright^{1}_{\l}, \vartriangleleft^{1}_{\l})$ and $(A, \vartriangleright^{2}_{\l}, \vartriangleleft^{2}_{\l})$ are compatible L-dendriform conformal algebras, where for all  $a, b \in A$, 
\begin{eqnarray*}
 &&a \vartriangleright^{1}_{\l} b = \mathcal{R}(a)\circ_{\l}b, ~\quad\quad a \vartriangleleft^{1}_{\l} b = -b\circ_{-\partial-\l} \mathcal{R}(a),\\
&&a \vartriangleright^{2}_{\l} b = (N\mathcal{R})(a)\circ_{\l}b, \quad a \vartriangleleft^{2}_{\l} b = -b\circ_{-\partial-\l} (N\mathcal{R})(a).   
\end{eqnarray*}

If in addition, $\mathcal R$ is invertible, then $(A, \vartriangleright^{1}_{\l}, \vartriangleleft^{1}_{\l})$ and $(A, \vartriangleright^{2}_{\l}, \vartriangleleft^{2}_{\l})$ are also compatible L-dendriform conformal algebras, with the operations defined for all $a, b \in A$ by    
\begin{eqnarray*}
 &&a \vartriangleright^{1}_{\l} b = \mathcal{R}(a\circ_{\l}\mathcal{R}^{-1}(b)), \quad \quad \quad ~~ a \vartriangleleft^{1}_{\l} b = -\mathcal{R}(\mathcal{R}^{-1}(b)\circ_{-\partial-\l} a),\\
&&a \vartriangleright^{2}_{\l} b = (N\mathcal{R})(a\circ_{\l}(N\mathcal{R})^{-1}(b)), \quad a \vartriangleleft^{2}_{\l} b = -(N\mathcal{R})((N\mathcal{R})^{-1}(b)\circ_{-\partial-\l} a).  
\end{eqnarray*}
    
\end{ex}

\end{document}